\documentclass{amsart}

\usepackage[utf8]{inputenc}
\usepackage[T1]{fontenc}
\usepackage{srcltx}
\usepackage{amssymb}
\usepackage{amsmath}
\usepackage{amsfonts}
\usepackage{latexsym}
\usepackage{amsthm}
\usepackage{enumerate}

\newtheorem{theo}{Theorem}[section]
\newtheorem{lemm}[theo]{Lemma}

\newtheorem{defi}[theo]{Definition}
\newtheorem{propo}[theo]{Proposition}
\newtheorem{coro}[theo]{Corollary}

\def\Z{\mathbb{Z}}
\def\N{\mathbb{N}}
\def\F{\mathbb{F}}

\def\Q{\mathbb{Q}}
\def\R{\mathbb{R}}
\def\C{\mathbb{C}}

\newcommand{\Ld}[2][D]{\ensuremath{\mathcal{L}(#2\mathcal{#1})}} 
\newcommand{\D}[1][D]{\ensuremath{\mathcal{#1}}} 

\begin{document}

\title[]{On the tensor rank
of multiplication in finite extensions of finite fields}
\author{S. Ballet}
\address{Institut de Math\'{e}matiques
de Luminy\\ case 930, F13288 Marseille cedex 9\\ France}
\email{stephane.ballet@univmed.fr}
\author{J. Chaumine}
\address{Laboratoire G\'eom\'etrie Alg\'ebrique et 
Applications \`a la Th\'eorie de l'Information\\ 
Universit\'e de la Polyn\'esie Fran\c{c}aise,\\ 
B.P. 6570, 98702 Faa'a, Tahiti\\ France}
\email{jean.chaumine@upf.pf}
\author{J. Pieltant}
\address{Institut de Math\'{e}matiques
de Luminy\\ case 930, F13288 Marseille cedex 9\\ France}
\email{julia.pieltant@univmed.fr}
\author{R. Rolland}
\address{Institut de Math\'{e}matiques
de Luminy\\ case 930, F13288 Marseille cedex 9\\ France}
\email{robert.rolland@acrypta.fr}

\date{\today}

\begin{abstract}
In this paper, we give a survey of the known results concerning the tensor rank of the multiplication 
in finite fields and we establish new asymptotical and not asymptotical upper bounds about it. 
\end{abstract}

\maketitle


\section{Introduction}
Several objects constitute the aim of this paper. 
First, it is a question of introducing the problem of the tensor rank of the multiplication 
in finite fields and of giving  a statement of the results obtained in this
part of algebraic complexity theory for which the best general reference is \cite{buclsh}.  
In particular, one of the aims of this paper  is to list exhaustively the few published mistaken 
statements and to explain them.  In the second part, we repair and clarify certain of these statements.
Last but not least, we improve several known results. 
In this section we introduce the problem, we set up notation and terminology
and we present the organization of this paper as well as the new obtained results.

\subsection{The bilinear complexity of the multiplication}

Let $\F_q$ be a finite field with $q=p^r$ elements where $p$ is a prime number. 
Let $\F_{q^n}$ be a degree $n$ extension of $\F_q$.
The multiplication $m$ in the finite field $\F_{q^n}$ is a bilinear map 
from $\F_{q^n} \times \F_{q^n}$ into $\F_{q^n}$, 
thus it corresponds to a linear map $M$ from the tensor product 
$\F_{q^n} \bigotimes \F_{q^n}$ into $\F_{q^n}$.
One can also represent $M$ by a tensor 
$t_M \in \F_{q^n}^*\bigotimes \F_{q^n}^* \bigotimes \F_{q^n}$ 
where $\F_{q^n}^*$ denotes the algebraic dual of $\F_{q^n}$.
Each decomposition
\begin{equation}\label{algo}
t_M=\sum_{i=1}^{k}a^*_i\otimes b^*_i\otimes c_i
\end{equation}
of the tensor $t_M$, where $a^*_i, b^*_i \in \F_{q^n}^*$ and
$c_i \in \F_{q^n}$, brings forth a multiplication algorithm
$$x.y=t_M(x\otimes y)=\sum_{i=1}^{k}a^*_i(x)\otimes b^*_i(x)\otimes c_i.$$

The bilinear complexity of the multiplication in $\F_{q^n}$ over $\F_q$, 
denoted by $\mu_{q}(n)$,
is the minimum number of summands in the decomposition (\ref{algo}).
Alternatively, we can say that the bilinear complexity of the multiplication 
is the rank of the tensor $t_M$ (cf. \cite{shtsvl}, \cite{ball2}).

\subsection{Organization of the paper}

In Section 2, we present the classical results via the approach using the multiplication by polynomial interpolation.
In section 3, we give an historical record of results resulting from the pioneer works due to 
D.V. and G.V. Chudnovky \cite{chch} and later Shparlinski, Tsfasman and Vladut in \cite{shtsvl}. 
In particular in Subsection 3.1, we present the original algorithm as well as 
the most successful version of the algorithm of Chudnovsky type at the present time. 
This modern approach uses the interpolation over algebraic curves defined over 
finite fields. This approach, which we recount the first success 
as well as   the rocks on which the pionners came to grief, enables to end at a first complete proof of the linearity of 
the bilinear complexity of multiplication \cite{ball1}. Then, in Subsection 3.2, 
we recall  the known results about the bilinear complexity $\mu_q (n)$.
Finally, in Section 4, we give new results for $\mu_q(n)$. 
More precisely, we obtain new upper bounds for $\mu_q(n)$ 
as well as new asymptotical upper bounds.

\section{Old classical results}
Let 
$$P(u)=\sum_{i=0}^{n}a_iu^i$$
be a  monic irreducible polynomial of degree $n$ with
coefficients in a field $F$.
Let
$$R(u)=\sum_{i=0}^{n-1}x_iu^i$$
and
$$S(u)=\sum_{i=0}^{n-1}y_iu^i$$
be two polynomials of degree $\leq n-1$ where
the coefficients $x_i$ and $y_i$ are indeterminates.

Fiduccia and Zalcstein (cf. \cite{fiza}, \cite{buclsh}  p.367 prop. 14.47)
have studied the general problem of computing
the coefficients of the product $R(u) \times S(u)$
and they have shown that at least $2n-1$ multiplications
are needed.
When the field $F$ is infinite, an algorithm reaching exactly
this bound was previously given by Toom in \cite{toom}.
Winograd described in \cite{wino2} all the algorithms reaching
the bound $2n-1$.
Moreover, Winograd proved in \cite{wino3} that up to some
transformations every algorithm for computing the coefficients
of $R(u) \times S(u) \mod P(u)$ which is of bilinear complexity $2n-1$,
necessarily computes the coefficients of $R(u) \times S(u)$, and consequently
uses one of the algorithms described in \cite{wino2}.
These algorithms use interpolation technics and cannot be performed
if the cardinality of the field $F$ is $<2n-2$. In conclusion
we have the following result:

\begin{theo}\label{old}
If the cardinality of $F$ is $<2n-2$, every algorithm computing
the coefficients of  $R(u) \times S(u) \mod P(u)$ has a bilinear
complexity $>2n-1$.
\end{theo}

Applying the results of Winograd and De Groote \cite{groo} and Theorem \ref{old} to
the multiplication in a finite extension $\F_{q^n}$ of a finite field
$\F_q$ we obtain:

\begin{theo}\label{thm_wdg}
The bilinear complexity $\mu_q(n)$ of the multiplication
in the finite field $\F_{q^n}$ over $\F_q$ verifies
$$\mu_q(n) \geq 2n-1,$$
with equality holding if and only if
$$n \leq \frac{q}{2}+1.$$
\end{theo}

This result does not give any estimate of an upper bound for 
$\mu_q(n)$, when $n$ is large.
In \cite{lesewi}, Lempel, Seroussi and Winograd proved
that $\mu_q(n)$ has a quasi-linear upper bound.
More precisely:

\begin{theo}
The bilinear complexity of the multiplication
in the finite field $\F_{q^n}$ over $\F_q$ verifies:
$$\mu_q(n) \leq f_q(n)n,$$
where $f_q(n)$ is a very slowly growing function, namely
$$f_q(n)=O(\underbrace{\log_q\log_q \cdots \log_q}_{k~times}(n))$$
for any $k \geq 1$.
\end{theo}

Furthermore, extending and using more efficiently the technique
developed in \cite{bska1}, Bshouty and Kaminski
showed that
$$\mu_q(n) \geq 3n-o(n)$$
for $q \geq 3.$
The proof of the above lower bound on the complexity of straight-line
algorithms for polynomial multiplication is based on the analysis
of Hankel matrices representing bilinear forms defined by linear
combinations of the coefficients of the polynomial product.

\section{The modern approach via algebraic curves}

We have seen in the previous section
that if the number of points of the ground field is too low, we cannot perform
the multiplication by the Winograd interpolation method. 
D.V. and G.V. Chudnowsky
have designed in \cite{chch} an algorithm 
where the interpolation is done on points of an 
algebraic
curve over the groundfield with a sufficient number of rational points. 
Using this algorithm, D.V. and G.V. Chudnovsky claimed that 
the bilinear complexity of
the multiplication in finite extensions of a finite field is asymptotically linear but later Shparlinski, Tsfasman and Vladut in \cite{shtsvl}  
noted that they only proved that the quantity $m_q=\liminf_{k \rightarrow \infty}\frac{\mu_q(k)}{k}$ is bounded 
which do not enable to prove the linearity. To prove the linearity, it is also necessary to prove that 
$M_q= \limsup_{k \rightarrow \infty}\frac{\mu_q(k)}{k}$ is bounded which is the main aim of their paper. 
However,  I. Cascudo, R. Cramer and C. Xing recently detected a mistake in the proof of Shparlinski, Tsfasman and Vladut. 
Unfortunately, this mistake that we will explain in details in this section, also had an effect on their improved estimations of $m_q$. 
After the above pioneer research, S. Ballet obtained in \cite{ball1} the first upper bounds uniformly with respect to $q$ for $\mu_q(n)$. 
These bounds not being affected by the same mistake enable at the same time to prove the linearity of the bilinear complexity of
the multiplication in finite extensions of a finite field. Then,  S. Ballet and al. obtained several improvements which will be recalled at the end of this section. These different improvements are based on the following main ideas: the use of towers of algebraic functions fields \cite{ball1} \cite{ball3}, the descent of their definition field \cite{baro1} \cite{balbro}, the use of places of higher degree \cite{baro1}  \cite{ceoz} as well as the use of local expansion \cite{arna1} \cite{ceoz}.

\subsection{Linearity of the bilinear complexity of the multiplication}

\subsubsection{The D.V. Chudnovsky and G.V. Chudnovsky algorithm}

In this section, we recall the brilliant idea of  D.V. Chudnovsky and G.V. Chudnovsky and give their main result.
First, we present  the original algorithm of D.V. Chudnovsky and G.V. Chudnovsky, which was established in 1987 in \cite{chch}.

\begin{theo}\label{chudchud}
Let 

\vspace{.1em}

$\bullet$ $F/\F_q$ be an algebraic function field,

\vspace{.1em}

$\bullet$ $Q$ be a degree $n$ place of $F/\F_q$, 

\vspace{.1em}

$\bullet$ ${\mathcal D}$ be a divisor of $F/\F_q$, 

\vspace{.1em}

$\bullet$ ${\mathcal P}=\{P_1,...,P_{N}\}$ 
be a set of places of degree $1$.

\vspace{.1em}

We suppose that $Q$, $P_1, \cdots, P_N$ are not in the support of ${\mathcal D}$
and that:

\vspace{.1em}

a) The evaluation map 
$$Ev_Q: {\mathcal L}({\mathcal D}) \rightarrow \F_{q^n}\simeq F_Q$$
is onto (where $F_Q$ is the residue class field of $Q$),

\vspace{.1em}

b) the application 

$$Ev_{{\mathcal P}}:
\left \{
\begin{array}{lll}
  {\mathcal L}(2{\mathcal D}) & \rightarrow & \F_{q}^{N}  \cr
   f & \mapsto & (f(P_1),...,f(P_{N})) \cr
\end{array}
\right .
$$

is injective.

\vspace{.1em}

Then  $$\mu_q(n)\leq N.$$
\end{theo}

\medskip

As pointed in \cite{shtsvl}, using this algorithm with a suitable sequence of algebraic curves defined over a finite field $\F_q$,
D.V. Chudnovsky and G.V. Chudnovsky only proved the following result:

\begin{theo}\label{chudmq}
Let $q$ be a square $\geq 25$. Then

$$\liminf \frac{\mu_q(n)}{n}\leq 2\left(1+\frac{1}{\sqrt{q}-3}\right).$$
\end{theo}

\subsubsection{Asymptotic bounds}\label{mM}
As seen previously, Shparlinski, Tsfasman, Vladut have given in \cite{shtsvl} many interesting
remarks on the algorithm of D.V. and G.V. Chudnovsky and the bilinear complexity.
In particular, they have considered asymptotic bounds
for the bilinear complexity in order to prove the asymptotic linearity of this complexity 
from the algorithm of D.V. and G.V. Chudnovsky.
Following these authors, let us define 
$$M_q= \limsup_{k \rightarrow \infty}\frac{\mu_q(k)}{k}$$
and
$$m_q=\liminf_{k \rightarrow \infty}\frac{\mu_q(k)}{k}.$$

It is not at all obvious that either of these values is finite but anyway the bilinear complexity of multiplication 
can be considered as asymptotically linear in the degree of extension if and only if the quantity $M_q$ is finite. 
First, let us recall a very useful Lemma due to  D.V. and G.V. Chudnovsky \cite{chch} and 
Shparlinski, Tsfasman, Vladut \cite[Lemma 1.2 and Corollary 1.3]{shtsvl}.

\begin{lemm}\label{lemasyMqmq}
For any prime power $q$ and for all the positive integers $n$ and $m$, we have 
$$ \mu_{q}(m)\leq\mu_q(mn)\leq \mu_q(n).\mu_{q^n}(m)$$
$$ m_{q}\leq m_{q^n}.\mu_{q}(n)/n$$
$$ M_{q}\leq M_{q^n}.\mu_{q}(n).$$
\end{lemm}

Now, let us summarize the known estimates concerning these quantities, namely the lower bound of $m_2$ 
obtained by R. Brockett, M. Brown and D. Dobkin in \cite{brttdo} \cite{brdo} 
and the lower bound of $m_q$ for $q>2$ given by Shparlinski, Tsfasman and Vladut in \cite{shtsvl}.

\begin{propo}

$$m_2\geq 3.52$$ and $$m_q\geq 2\left( 1+\frac{1}{q-1}\right) \hbox{ for any }q>2.$$

\end{propo}

\medskip

Note that all the upper bounds of $M_q$ and $m_q$ for any $q$ given by Shparlinski, 
Tsfasman and Vladut in \cite{shtsvl} are not proved.
Indeed, in \cite{shtsvl}, they claim that for any $q$ (in particular for $q=2$), $m_q$ and overall $M_q$ are finite but I. Cascudo, R. Cramer and C. Xing recently communicated us the existence of a gap in the proof established by I. Shparlinsky, M. Tsfasman and S. Vladut: \textit{"the mistake in \cite{shtsvl}  from 1992 is in the proof of their Lemma 3.3, page 161, the paragraph following formulas about the degrees of the divisor. It reads: "}\textsl{Thus the  number of linear equivalence classes of degree a for which either Condition $\alpha$ or Condition $\beta$ fails is at most $D_{b'} + D_b$.}\textit{" This is incorrect; $D_b$ should be multiplied by the torsion. Hence the proof of their asympotic bound is incorrect."} \\
Let us explain this gap in next section.

\subsubsection{Gap in the proof of the asymptotic linearity}\label{gap}

We settle the following elements 
\begin{enumerate}
\item a place of degree $n$ denoted by $Q$;
\item $2n+g-1$ places of degree $1$ : $P_1,\cdots,P_{2n+g-1}$.
\end{enumerate}

We look for a divisor $D$ such that:
\begin{enumerate}
\item $\deg(D)=n+g-1$;
\item $\dim({\mathcal L}(D-Q))=0$;
\item $\dim({\mathcal L}(2D - (P_1+P_2+\cdots+ P_{2n+g-1})))=0$.
\end{enumerate}

The results concerning $M_q$ et $m_q$ obtained in the paper \cite{tsvl} depend
on the existence of such a divisor $D$.

\medskip

Let us remark that these conditions only depend on the class of a divisor 
(the dimension of a divisor, the degree of a divisor are invariant in a same class). 
Consequently, we can work on classes and show the existence 
of a class $[D]$ which answers the question.

Let $J_{n+g-1}$ be the set of classes of degree $n+g-1$ divisors.
We know from F. K. Schmidt Theorem that there exists
a divisor $D_0$ of degree $n+g-1$.
The application
$\psi_{n+g-1}$ from $J_{n+g-1}$ into the Jacobian $J_0$
defined by
$$\psi_{n+g-1}([D])=[D-D_0]$$
is a bijection from $J_{n+g-1}$ into $J_0$.
All the sets $J_k$ have the same number $h$
of elements ($h$ is called the number of classes).

\medskip

Let $u$ be the application from $J_{n+g-1}$ into $J_{g-1}$
defined by $u([D])=[D-Q]$.
This application is bijective. 
Thus if we set 
$$H_{n+g-1}= \{ [D] \in J_{n+g-1} ~|~ \dim([D-Q])=0\},$$
and
$$K_{g-1}=\{ [{\Delta}] \in J_{g-1} ~|~ \dim([{\Delta}])=0\},$$
we have
$$K_{g-1}=u( H_{n+g-1}),$$
and then
$$\# H_{n+g-1} = \# K_{g-1}.$$

Let us note that if $[{\Delta}]$ is an element of $J_{g-1}$
which is in the complementary of $K_{g-1}$ namely
$\dim([{\Delta}])>0$, then there exists in the class $[{\Delta}]$
at least an effective divisor (there exists a $x$ such that 
${\Delta}+(x) \geq 0$).
Moreover effective divisors in different classes are different.
So the complementary of $K_{g-1}$ in $J_{g-1}$ has a cardinality
$\leq A_{g-1}$ where $A_{g-1}$ is the number of effective divisors of degree $g-1$.
Then the cardinality of $K_{g-1}$ verifies
the inequality
$$\# H_{n+g-1} = \# K_{g-1} \geq h-A_{g-1}.$$

\medskip

Let us remark that classes which belong to $H_{n+g-1}$ are
the only ones which can solve our problem.
But they also have to verify the additional condition 
$$\dim({\mathcal L}(2D - (P_1+P_2+\cdots+ P_{2n+g-1})))=0.$$
We would like to use a combinatorial proof as for the first condition.

So we have to consider the application $v$ from $H_{n+g-1}$
to $J_{g-1}$ defined by
$$v([D])= [2D - (P_1+P_2+\cdots+ P_{2n+g-1})].$$
Unfortunately the application $[D]\mapsto [2D]$ is not necessarily injective.
This is related to $2$-torsion points of the Jacobian.
The fact that the application $v$ is not injective does not allow us to 
conclude that there exists an image "big"  enough and use a combinatorial
argument like in the first part.

\subsection{Known results about the bilinear complexity $\mu_{q}(n)$}

\subsubsection{Extensions of the Chudnovsky algorithm}

In order to obtain good estimates for the bilinear complexity, S. Ballet has given in \cite{ball1} some easy to verify conditions allowing the use of the D.V. and G.V. Chudnovsky  algorithm. Then S. Ballet and R. Rolland have generalized in \cite{baro1} the algorithm using places of degree $1$ and~$2$.

Let us present the last version of this algorithm, which is a generalization of the algorithm of type Chudnovsky introduced by N. Arnaud in \cite{arna1} and M. Cenk and F. \"Ozbudak in \cite{ceoz}. This generalization uses several coefficients in the local expansion at each place $P_i$ instead of just the first one. Due to the way to obtain the local expansion of a product from the local expansion of each term, the bound for the bilinear complexity involves the complexity notion $\widehat{M_q}(u)$ introduced by M. Cenk and F. \"Ozbudak in \cite{ceoz} and defined as follows:
\begin{defi}
We denote by $\widehat{M_q}(u)$ the minimum number of multiplications needed in $\F_q$ in order to obtain coefficients of the product of two arbitrary $u$-term polynomials modulo $x^u$ in $\F_q[x]$.
\end{defi}
For instance, we know that for all prime powers $q$, we have $\widehat{M_q}(2) \leq 3$ by \cite{ceoz2}.\\

Now we introduce the generalized algorithm of type Chudnovsky described in~\cite{ceoz}.

\begin{theo} \label{theo_evalder}
Let \\
\vspace{.1em}
$\bullet$ $q$ be a prime power,\\
\vspace{.1em}
$\bullet$ $F/\F_q$ be an algebraic function field,\\
\vspace{.1em}
$\bullet$ $Q$ be a degree $n$ place of $F/\F_q$,\\
\vspace{.1em}
$\bullet$ ${\mathcal D}$ be a divisor of $F/\F_q$,\\
\vspace{.1em}
$\bullet$ ${\mathcal P}=\{P_1,\ldots,P_N\}$ be a set of $N$ places of arbitrary degree,\\
\vspace{.1em}
$\bullet$ $u_1,\ldots,u_N$ be positive integers.\\
We suppose that $Q$ and all the places in $\mathcal P$ are not in the support of ${\mathcal D}$ and that:
\begin{enumerate}[a)]
	\item the map
	$$
	Ev_Q:  \left \{
	\begin{array}{ccl}
	\Ld{} & \rightarrow & \F_{q^n}\simeq F_Q\\
	f & \longmapsto & f(Q)
	\end{array} \right.
	$$ 
	is onto,
	\item the map
	$$
	Ev_{\mathcal P} :  \left \{
	\begin{array}{ccl}
	\Ld{2} & \longrightarrow & \left(\F_{q^{\deg P_1}}\right)^{u_1} \times \left(\F_{q^{\deg P_2}}\right)^{u_2} \times \cdots \times \left(\F_{q^{\deg P_N}}\right)^{u_N} \\
	f & \longmapsto & \big(\varphi_1(f), \varphi_2(f), \ldots, \varphi_N(f)\big)
	\end{array} \right.
	$$
	is injective, where the application $\varphi_i$ is defined by
	$$
	\varphi_i : \left \{
	\begin{array}{ccl}
	\Ld{2} & \longrightarrow & \left(\F_{q^{\deg P_i}}\right)^{u_i} \\
          f & \longmapsto & \left(f(P_i), f'(P_i), \ldots, f^{(u_i-1)}(P_i)\right)
	\end{array} \right.
 	$$
	with $f = f(P_i) + f'(P_i)t_i + f''(P_i)t_i^2+ \ldots + f^{(k)}(P_i)t_i^k + \ldots $, 
the local expansion at $P_i$ of $f$ in ${\Ld{2}}$, with respect to the local parameter~$t_i$. 
Note that we set ${f^{(0)} =f}$.
\end{enumerate}
Then 
$$
\mu_q(n) \leq \displaystyle \sum_{i=1}^N \mu_q(\deg P_i) \widehat{M}_{q^{\deg P_i}}(u_i).
$$
\end{theo}

Let us remark that the algorithm given in \cite{chch} by D.V. and G.V. Chudnovsky is the 
case $\deg P_i=1$ and $u_i=1$ for $i=1, \ldots, N$. 
The first generalization introduced by S.Ballet and R. Rolland in \cite{baro1} 
concerns the case $\deg P_i=1 \hbox{ or }2$ and $u_i=1$ for $i=1, \ldots, N$. 
Next, the generalization introduced by N. Arnaud in \cite{arna1} concerns 
the case $\deg P_i=1 \hbox{ or }2$ and $u_i=1\hbox{ or }2$  for $i=1, \ldots, N$. 
However, note that the work of N. Arnaud has never been published and contains 
few mistakes (mentioned below) which will be repared in this paper. 
Finally, the last generalization introduced by M. Cenk and F. 
\"Ozbudak in \cite{ceoz} is useful: it allows us to use certain places of arbitrary degree many times, 
thus less places of fixed degree are necessary to get the injectivity of $Ev_\mathcal{P}$.

In particular, we have the following result, obtained by N. Arnaud in \cite{arna1}. 

\begin{coro} \label{theo_deg12evalder}
Let \\
\vspace{.1em}
$\bullet$ $q$ be a prime power,\\
\vspace{.1em}
$\bullet$ $F/\F_q$ be an algebraic function field,\\
\vspace{.1em}
$\bullet$ $Q$ be a degree $n$ place of $F/\F_q$,\\ 
\vspace{.1em}
$\bullet$ $\D$ be a divisor of $F/\F_q$, \\
\vspace{.1em}
$\bullet$ ${\mathcal P}=\{P_1,\ldots,P_{N_1},P_{N_1+1},\ldots,P_{N_1+N_2}\}$ be a set of 
$N_1$ places of degree\\ \indent one and $N_2$ places of degree two,\\
\vspace{.1em}
$\bullet$ ${0 \leq l_1 \leq N_1}$ and ${0 \leq l_2 \leq N_2}$ be two integers.\\
We suppose that $Q$ and all the places in $\mathcal P$ are not in the support of $\D$ and that:
\begin{enumerate}[a)]
	\item the map
	$$
	Ev_Q: \Ld{} \rightarrow \F_{q^n}\simeq F_Q$$
	is onto,
	\item the map
	 $$
	 Ev_{\mathcal P}: \left \{
	\begin{array}{ccl}
   	\Ld{2} & \rightarrow & \F_{q}^{N_1} \times \F_{q}^{l_1}\times \F_{q^2}^{N_2} \times \F_{q^2}^{l_2}  \\
    	f & \mapsto & \big(f(P_1),\ldots,f(P_{N_1}),f'(P_1),\ldots,f'(P_{l_1}),\\
	  &  & \ f(P_{N_1+1}),\ldots,f(P_{N_1+N_2}),f'(P_{N_1+1}),\ldots,f'(P_{N_1+l_2})\big)
	\end{array} \right .
	$$
is injective.
\end{enumerate}
Then  
$$
\mu_q(n)\leq  N_1 + 2l_1 + 3N_2 + 6l_2.
$$
\end{coro}

Moreover, from the last corollary applied on Garcia-Stichtenoth towers, N. Arnaud obtained in \cite{arna1} the two following bounds:

\begin{theo}\label{bornes_arnaud1}
Let ${q=p^r}$ be a prime power.
	\begin{equation*}
	\begin{array}{l}
	\mbox{(i) \ If $q\geq4$, then  }\mu_{q^2}(n) \leq 2 \left(1 + \frac{p}{q-3 + (p-1)\left(1- \frac{1}{q+1}\right)} \right)n,\\
	\\
	\mbox{(ii) \ If $q\geq16$, then }\mu_{q}(n) \leq 3 \left(1 + \frac{2p}{q-3 + 2(p-1)\left(1- \frac{1}{q+1}\right)} \right)n.
	\end{array}
	\end{equation*}
\end{theo}

We will give a proof of Bound (i) together with an improvement of Bound (ii) in Section \ref{sectbornesarnaud}. In that section, we will also prove two revised bounds for $\mu_{p^2}(n)$ and $\mu_p(n)$ given by Arnaud in \cite{arna1}. Indeed, Arnaud gives the two following bounds with no detailed calculation:
	\begin{equation*}
	\begin{array}{l}
	\mbox{(iii) \ If $p\geq5$ is a prime, then  }\mu_{p^2}(n) \leq 2 \left(1 + \frac{2}{p-2} \right)n,\\
	\\
	\mbox{(iv) \ If $p\geq5$ is a prime, then } \mu_p(n)\leq 3 \left(1+ \frac{4}{p-1}\right) n.
	\end{array}
	\end{equation*}
In fact, one can check that the denominators $p-1$ and $p-2$ are slightly overestimated under Arnaud's hypotheses.

From the 
results of \cite{ball1} and the previous
algorithm, we obtain (cf. \cite{ball1}, \cite{baro1}):

\begin{theo} \label{theoprinc}
Let $q$ be a prime power and let $n$ be an integer $>1$. 
Let $F/\F_q$ be an algebraic function field of genus $g$ 
and $N_k$ the number of places of degree $k$ in $F/\F_q$.
If $F/\F _q$ is such that $2g+1 \leq q^{\frac{n-1}{2}}(q^{\frac{1}{2}}-1)$ then:
\begin{enumerate}[1)]
	\item if $N_1 > 2n+2g-2$, then $$ \mu_q(n) \leq 2n+g-1,$$
	\item if there exists a non-special divisor of degree $g-1$ 
and $N_1+2N_2>2n+2g-2$, then $$\mu_q(n)\leq 3n+3g,$$
	\item if $N_1+2N_2>2n+4g-2$, then $$\mu_q(n)\leq 3n+6g.$$
\end{enumerate}
\end{theo}

\subsubsection{Known upper bounds for $\mu_{q}(n)$}

From good towers of algebraic functions fields satisfying Theorem \ref{theoprinc}, it 
was proved in 
\cite{ball1}, \cite{ball3}, \cite{baro1}, \cite{balbro}, \cite{ball4} and \cite{bach}:

\begin{theo}
Let $q=p^r$ a power of the prime $p$.
The bilinear complexity $\mu_q(n)$ of multiplication in any finite field $\F_{q^n}$  
is linear 
with respect to the extension degree, more precisely:
$$\mu_q(n) \leq C_q n$$ where $C_q$ is the constant defined by: 
$$
C_q=
\left \{
\begin{array}{lll}
\hbox{if } q=2 & \hbox{then} \quad 22 &  \hbox{\cite{bapi} and \cite{ceoz}} \cr \cr
\hbox{else if } q=3 &  \hbox{then} \quad 27 & \hbox{\cite{ball1}} \cr \cr
\hbox{else if } q=p \geq 5 &  \hbox{then} \quad  3\left(1+ \frac{4}{q-3}\right) & 
                   \hbox{\cite{bach}} \cr \cr
\hbox{else if } q=p^2 \geq 25 & \hbox{then} \quad  2\left(1+\frac{2}{p-3}\right) & 
                   \hbox{\cite{bach}} \cr  \cr
\hbox{else if } q=p^{2k} \geq 16 & \hbox{then} \quad 2\left(1+\frac{p}{q-3 + (p-1)\left(1- \frac{1}{q+1}\right)}\right) & 
                   \hbox{\cite{arna1}} \cr \cr
\hbox{else if } q \geq 4 & \hbox{then} \quad 6\left(1+\frac{p}{q-3}\right) & \hbox{\cite{ball3}}\cr \cr
\hbox{else if } q \geq 16 & \hbox{then} \quad 3\left(1+\frac{2p}{q-3 + 2(p-1)\left(1- \frac{1}{q+1}\right)}\right) & \hbox{\cite{arna1}}.
\end{array}
\right .
$$
\end{theo}

Note that the new estimate for the constant $C_2$ comes from two recent improvements. First, one knows from Table 1 in \cite{ceoz} that $\mu_2(n) \leq 22n$  for ${2 \leq n \leq 7}$ since $\mu_2(n) \leq 22$ for such integers $n$. Moreover, applying the bound ${\mu_2(n)\leq \frac{477}{26}n+\frac{45}{2}}$ obtained in \cite{bapi}, one gets  ${\mu_2(n)\leq \left(\frac{477}{26}+\frac{45}{2\times8}\right)n \leq 22n}$ for ${n\geq 8}$. Note also that the upper bounds obtained in \cite{ball5} and \cite{ball6} 
are obtained by using the mistaken statements of  I. Shparlinsky, M. Tsfasman and S. Vladut \cite{shtsvl} mentioned 
in the above section \ref{gap}. Consequently, these bounds are not proved and unfortunatly they can not be repaired easily. 
However, certain not yet published results recently due to H. Randriambololona concerning the geometry of Riemann-Roch spaces 
might enable to repair them in certain cases.

\subsubsection{Some exact values for the bilinear complexity}

Applying the D.V. and G.V. Chudnovsky algorithm with
well fitted elliptic curves, Shokrollahi has shown in \cite{shok}
that:

\begin{theo}\label{thm_shokr}
The bilinear complexity $\mu_q(n)$ of 
the multiplication in the finite extention $\F_{q^n}$ of the
finite field $\F_q$ is equal to $2n$ for 
\begin{equation}\label{ine}
\frac{1}{2}q +1< n < \frac{1}{2}(q+1+{\epsilon (q) })
\end{equation} 
where $\epsilon$ is the function defined by:
$$
\epsilon (q)=
\left \{
\begin{array}{l}
 \hbox{the greatest integer} \le 2{\sqrt q} \hbox{ prime to q,} 
 \quad \hbox{if}~ \hbox{q is not a perfect square} \cr
 2{\sqrt q,} \quad \hbox{if}~ \hbox{q is a perfect square.} \cr
\end{array}
\right .
$$
\end{theo}

We still do not know if the converse is true. More precisely the question is:
suppose that $\mu_q(n)=2n$, are the inequalities (\ref{ine}) true?

However, for computational use, it is helpful to keep in mind some particular exact values for ${\mu_q(n)}$, such as ${\mu_q(2)=3}$ for any prime power $q$, ${\mu_2(4)=9}$, ${\mu_4(4)=\mu_5(4)=8}$ or ${\mu_2(2^6)=15}$ \cite{chch}.

\section{New results for $\mu_q(n)$}

\subsection{Towers of algebraic function fields}\label{sectdeftowers}

In this section, we introduce some towers of algebraic function fields.
Theorem \ref{theoprinc} applied to the algebraic function fields
of these towers gives us bounds for the bilinear complexity. 
A given curve cannot permit to multiply in every extension of $\F_q$,
just for $n$ lower than some value. With a tower of function fields
we can adapt the curve to the degree of the extension. 
The important point to note here is that in order
to obtain a well adapted curve it will be desirable
to have a tower for which the quotients of two consecutive genus
are as small as possible, namely a dense tower. 

\medskip

For any algebraic function field $F/\F_q$ defined over the finite field $\F_q$, 
we denote by $g(F/\F_q)$ the genus of $F/\F_q$
and by $N_k(F/\F_q)$ the number of places of degree $k$ in $F/\F_q$.

\subsubsection{Garcia-Stichtenoth tower of Artin-Schreier algebraic
function field extensions}

We present now a modified 
Garcia-Stichtenoth's tower (cf. \cite{gast}, \cite{ball3}, \cite{baro1}) 
having good properties. 
Let  us consider a finite field $\F_{q^2}$ with $q=p^r>3$ and $r$ an odd integer. 
Let us consider the Garcia-Stichtenoth's elementary abelian tower $T_1$ over $\F_{q^2}$ 
constructed in \cite{gast} and defined by the sequence $(F_0, F_1, F_2,\ldots)$ where 
$$F_{k+1}:=F_{k}(z_{k+1})$$ 
and $z_{k+1}$ 
satisfies the equation: 
$$z_{k+1}^q+z_{k+1}=x_k^{q+1}$$ 
with 
$$x_k:=z_k/x_{k-1} ~ in ~ F_k(for~k\geq1).$$
Moreover $F_0:=\F_{q^2}(x_0)$ is the rational function field over $\F_{q^2}$ 
and  $F_1$ the Hermitian function field over $\F_{q^2}$.
Let us denote by $g_k$ the genus of $F_k$, we recall the following \textsl{formulae}:
\begin{equation}\label{genregs}
g_k = \left \{ \begin{array}{ll}
		q^k+q^{k-1}-q^\frac{k+1}{2} - 2q^\frac{k-1}{2}+1 & \mbox{if } k \equiv 1 \mod 2,\\
		q^k+q^{k-1}-\frac{1}{2}q^{\frac{k}{2}+1} - \frac{3}{2}q^{\frac{k}{2}}-q^{\frac{k}{2}-1} +1& \mbox{if } k \equiv 0 \mod 2.
		\end{array} \right .
\end{equation}
Let us consider the completed Garcia-Stichtenoth tower 
$$T_2=F_{0,0}\subseteq F_{0,1}\subseteq \ldots \subseteq F_{0,r} \subseteq F_{1,0}\subseteq F_{1,1} \subseteq \ldots \subseteq F_{1,r} \ldots $$ 
considered in \cite{ball3} such that   
$F_k \subseteq F_{k,s} \subseteq F_{k+1}$ for any integer $s \in \{0,\ldots,r\}$,  
with $F_{k,0}=F_k$ and $F_{k,r}=F_{k+1}$. Recall that each extension $F_{k,s}/F_k$ is Galois of degree $p^s$ 
with full constant field $\F_{q^2}$. 
Now, we consider the tower studied in \cite{baro1}
$$T_3=G_{0,0} \subseteq G_{0,1} \subseteq \ldots \subseteq G_{0,r}\subseteq G_{1,0}\subseteq G_{1,1}\subseteq \ldots \subseteq G_{1,r}\ldots $$
defined over the constant field $\F_q$ and related to
the tower $T_2$ by
$$F_{k,s}=\F_{q^2}G_{k,s} \quad \mbox{for all $k$ and $s$,}$$
namely $\F_{k,s}/\F_{q^2}$ is the constant field extension of $G_{k,s}/\F_q$. Note that the tower $T_3$ is well defined 
by \cite{baro1} and \cite{balbro}. Moreover, we have the following result:

\begin{propo}\label{subfield}
Let ${q = p^r \geq4}$ be a prime power. For all integers $k \geq 1$ and ${s \in \{0, \ldots,r\}}$, there exists a step $F_{k,s}/\F_{q^2}$ (respectively $G_{k,s}/\F_q$) with genus $g_{k,s}$ and $N_{k,s}$ places of degree 1 in $F_{k,s}/\F_{q^2}$ (respectively  $N_{k,s}$ places of degree 1 and 2 in $G_{k,s}/\F_q$ with places of degree 2 being counted twice) such that:
\begin{enumerate}[(1)]
	\item $F_k \subseteq F_{k,s} \subseteq F_{k+1}$, where we set $F_{k,0}=F_k$ and $F_{k,r}=F_{k+1}$,\\(respectively $G_k \subseteq G_{k,s} \subseteq G_{k+1}$, where we set $G_{k,0}=G_k$ and $G_{k,r}=G_{k+1}$),
	\item $\big( g_k-1 \big)p^s +1 \leq g_{k,s} \leq \frac{g_{k+1}}{p^{r-s}} +1$,
	\item $N_{k,s} \geq (q^2-1)q^{k-1}p^s$.
\end{enumerate} 
\end{propo}

\subsubsection{Garcia-Stichtenoth tower of Kummer function field extensions}

In this section we present a Garcia-Stichtenoth's tower 
(cf. \cite{bach}) having good properties.
Let $\F_q$ be a finite field of characteristic $p\geq3$.
Let us consider the tower $T$ over $\F_q$ which is defined recursively 
by the following equation, studied in \cite{gast2}:

$$y^2=\frac{x^2+1}{2x}.$$

The tower $T/\F_q$ is represented by 
the sequence of function fields $(H_0, H_1, H_2, ...)$ 
where $H_n = \F_q(x_0, x_1, ..., x_n)$ and $x_{i+1}^2=(x_i^2+1)/2x_i$ 
holds for each $i\geq 0$.
Note that $H_0$ is the rational function field.
For any prime number $p \geq 3$, the tower $T/\F_{p^2}$ is 
asymptotically optimal over the field $\F_{p^2}$,
i.e. $T/\F_{p^2}$ reaches the Drinfeld-Vladut bound.
Moreover, for any integer $k$, $H_k/\F_{p^2}$ is the constant field extension of $H_k/\F_p$.

From \cite{bach}, we know that the genus $g(H_k)$ of the step $H_k$ is given by:
\begin{equation}\label{genregsr}
g(H_k) = \left \{ \begin{array}{ll}
		2^{k+1}-3\cdot 2^\frac{k}{2}+1 & \mbox{if } k \equiv 0 \mod 2,\\
		2^{k+1} -2\cdot 2^\frac{k+1}{2}+1& \mbox{if } k \equiv 1 \mod 2.
		\end{array} \right .
\end{equation}
and that the following bounds hold for the number of rational places  in $H_k$ over $\F_{p^2}$ and for  the number of places of degree 1 and 2 over $\F_p$:
\begin{equation}\label{nbratplgsr}
N_1(H_k/\F_{p^2}) \geq 2^{k+1}(p-1)
\end{equation}
and
\begin{equation}\label{nbpldeg12gsr}
N_1(H_k/\F_p) +2N_2(H_k/\F_p) \geq 2^{k+1}(p-1).
\end{equation}

From the existence of this tower, we can obtain the following proposition 
\cite{bach}:
\begin{propo} \label{propoexistcf2}
Let $p$ be a prime number $\geq 5$.
Then for any integer $n\geq {\frac{1}{2} (p+1+\epsilon(p))}$ 
where $\epsilon(p)$ is defined as in Theorem \ref{thm_shokr},
\begin{enumerate}[1)]
	\item there exists an algebraic function field $H_k/\F_{p^2}$ of genus $g(H_k/\F_{p^2})$
such that $2g(H_k/\F_{p^2})+1 \leq p^{{n-1}}(p-1)$
and $N_1(H_k/\F_{p^2})>2n+2g(H_k/\F_{p^2})-2$,
	\item there exists an algebraic function field $H_k/\F_{p}$ of genus $g(H_k/\F_p)$
such that $2g(H_k/\F_p)+1 \leq p^{\frac{n-1}{2}}(p^{\frac{1}{2}}-1)$
and $N_1(H_k/\F_p)+2N_2(H_k/\F_p)>2n+2g(H_k/\F_p)-2$
and containing a non-special divisor of degree $g(H_k/\F_p)-1$.
\end{enumerate}
\end{propo}

\subsection{Some preliminary results} \label{sectionusefull}
Here we establish some technical results about genus and number of places of each step of the towers $T_2/\F_{q^2}$, $T_3/\F_q$, $T/\F_{p^2}$ and $T/\F_p$ defined in Section \ref{sectdeftowers}. These results will allow us to determine a suitable step of the tower to apply the algorithm on. 
\subsubsection{About the Garcia-Stichtenoth's tower} In this section, $q:=p^r$ is a power of the prime $p$. 

\begin{lemm}\label{lemme_genre}
Let ${q>3}$. We have the following bounds for the genus of each step of the towers $T_2/\F_{q^2}$ and  $T_3/\F_q$:
\begin{enumerate}[i)]
	\item $g_k> q^k$ for all ${k\geq 4}$,
	\item $g_k \leq q^{k-1}(q+1) - \sqrt{q}q^\frac{k}{2}$,
	\item $g_{k,s} \leq q^{k-1}(q+1)p^s$ for all ${k\geq 0}$ and $s=0,\ldots,r$,
	\item $g_{k,s} \leq \frac{q^k(q+1)-q^\frac{k}{2}(q-1)}{p^{r-s}}$ for all $k\geq 2$ and $s=0,\ldots,r$.
\end{enumerate}
\end{lemm}

\begin{Proof}
\textit{i)} According to Formula (\ref{genregs}), we know that if ${k \equiv 1 \mod 2}$, then 
$$
g_k = q^k+q^{k-1}-q^\frac{k+1}{2} - 2q^\frac{k-1}{2}+1 = q^k+q^\frac{k-1}{2}(q^\frac{k-1}{2} - q - 2) +1.
$$
Since ${q>3}$ and ${k \geq 4}$, we have ${q^\frac{k-1}{2} - q - 2 >0}$, thus ${g_k>q^k}$.\\
Else if ${k \equiv 0 \mod 2}$, then 
$$
g_k = q^k+q^{k-1}-\frac{1}{2}q^{\frac{k}{2}+1} - \frac{3}{2}q^{\frac{k}{2}}-q^{\frac{k}{2}-1} +1 = q^k+q^{\frac{k}{2}-1}(q^\frac{k}{2}-\frac{1}{2}q^{2} - \frac{3}{2}q-1)+1.
$$
Since ${q>3}$ and ${k\geq 4}$, we have ${q^\frac{k}{2}-\frac{1}{2}q^{2} - \frac{3}{2}q-1>0}$, thus ${g_k>q^k}$.

\textit{ii)} It follows from Formula (\ref{genregs}) since for all $k\geq 1$ we have ${2q^\frac{k-1}{2} \geq 1}$ which works out for odd $k$ cases and ${\frac{3}{2}q^\frac{k}{2}+q^{\frac{k}{2}-1}\geq 1}$ which works out for even $k$ cases, since ${\frac{1}{2}q\geq \sqrt{q}}$.

\textit{iii)} If ${s=r}$, then according to Formula (\ref{genregs}), we have 
$$
g_{k,s} = g_{k+1}\leq q^{k+1}+q^{k} = q^{k-1}(q+1)p^s.
$$
Else, ${s<r}$ and Proposition \ref{subfield} says that ${g_{k,s} \leq \frac{g_{k+1}}{p^{r-s}}+1}$. Moreover, since ${q^\frac{k+2}{2}\geq q}$ and ${\frac{1}{2}q^{\frac{k+1}{2}+1}\geq q}$, we obtain ${g_{k+1}\leq q^{k+1} + q^k - q + 1}$ from Formula (\ref{genregs}). Thus, we get 
\begin{eqnarray*}
g_{k,s} & \leq & \frac{q^{k+1} + q^k - q + 1}{p^{r-s}} +1\\
	    &  = & q^{k-1}(q+1)p^s - p^s + p^{s-r} + 1\\
	    & \leq & q^{k-1}(q+1)p^s + p^{s-r}\\
	    & \leq & q^{k-1}(q+1)p^s \  \mbox{ since ${0 \leq p^{s-r} <1}$ and ${g_{k,s} \in \N}$}.
\end{eqnarray*}

\textit{iv)} It follows from ii) since Proposition \ref{subfield} gives ${g_{k,s} \leq \frac{g_{k+1}}{p^{r-s}}+1}$, so \linebreak[4]${g_{k,s} \leq  \frac{q^{k}(q+1) - \sqrt{q}q^\frac{k+1}{2}}{p^{r-s}} +1}$ which gives the result since ${p^{r-s} \leq q^\frac{k}{2}}$ for all ${k\geq2}$.
\qed\\
\end{Proof}

\begin{lemm}{\label{lemme_delta}}
Let $q>3$ and $k\geq4$. We set ${\Delta g_{k,s} := g_{k,s+1} - g_{k,s}}$ and ${D_{k,s}:=(p-1)p^sq^k}$ and denote ${M_{k,s} := N_1(F_{k,s}/\F_{q^2}) =  N_1(G_{k,s}/\F_q)+2N_2(G_{k,s}/\F_q)}$.
One has:
\begin{enumerate}[(i)]
	\item $\Delta g_{k,s} \geq D_{k,s}$,
	\item $M_{k,s} \geq D_{k,s}$.
\end{enumerate}
\end{lemm}

\begin{Proof}
(i) From  Hurwitz Genus Formula, one has ${g_{k,s+1}-1 \geq p(g_{k,s}-1)}$, so ${g_{k,s+1}-g_{k,s} \geq (p-1)(g_{k,s}-1)}$. Applying $s$ more times Hurwitz Genus Formula, we get ${g_{k,s+1}-g_{k,s} \geq (p-1)p^s\big(g(G_k)-1\big)}$. Thus ${g_{k,s+1}-g_{k,s} \geq (p-1)p^sq^k}$, from Lemma \ref{lemme_genre} i) since $q>3$ and $k\geq 4$.\\
\noindent(ii) According to Proposition \ref{subfield}, one has
\begin{eqnarray*}
	M_{k,s} & \geq & (q^2-1)q^{k-1}p^s \\
		    & = & (q+1)(q-1)q^{k-1}p^s \\
		    & \geq & (q-1)q^kp^s\\
		    & \geq & (p-1)q^kp^s \mbox{.}
\end{eqnarray*}
\qed
\end{Proof}

\begin{lemm}\label{lemme_bornesup}
Let ${M_{k,s} := N_1(F_{k,s}/\F_{q^2}) =  N_1(G_{k,s}/\F_q)+2N_2(G_{k,s}/\F_q)}$. For all ${k \geq 1}$ and ${s=0, \ldots, r}$, we have
$$
\sup \big \{ n \in \N \; | \;  2n \leq M_{k,s} -2g_{k,s} +1 \big \} \geq \frac{1}{2}(q+1)q^{k-1}p^s(q-3).
$$
\end{lemm}

\begin{Proof}
From Proposition \ref{subfield} and Lemma \ref{lemme_genre} iii), we get 
\begin{eqnarray*}
M_{k,s} -2g_{k,s} +1 & \geq & (q^2-1)q^{k-1}p^s - 2q^{k-1}(q+1)p^s +1 \\
				& = & (q+1)q^{k-1}p^s\big((q-1) -2\big) +1 \\
				& \geq & (q+1)q^{k-1}p^s(q-3)
\end{eqnarray*}
thus we have $\sup \big \{ n \in \N \; | \;  2n \leq M_{k,s} -2g_{k,s} +1 \big \} \geq \frac{1}{2}q^{k-1}p^s(q+1)(q-3)$.
\qed
\end{Proof}

\subsubsection{About the Garcia-Stichtenoth-R\"uck's tower} In this section, $p$ is an odd prime. We denote by $g_k$ the genus of the step $H_k$ and we fix $N_k:= N_1(H_k/\F_{p^2})=N_1(H_k/\F_p)+2N_2(H_k/\F_p)$. The following lemma is straightforward according to Formulae ~(\ref{genregsr}) and (\ref{nbpldeg12gsr}):
\begin{lemm}\label{lemme_genregsr}
These two bounds hold for the genus of each step of the towers $T/\F_{p^2}$ and $T/\F_p$:
\begin{enumerate}[i)]
	\item $g_k \leq 2^{k+1}-2\cdot2^\frac{k+1}{2}+1$,
	\item $g_k \leq 2^{k+1}$.
\end{enumerate}
\end{lemm}

\begin{lemm}\label{lemme_deltagsr}
For all $k\geq 0$, we set ${\Delta g_k := g_{k+1} - g_k}$. Then one has \linebreak[4]${N_k \geq \Delta g_k \geq 2^{k+1}- 2^\frac{k+1}{2}}$.
\end{lemm}

\begin{Proof} If $k$ is even then ${\Delta g_k = 2^{k+1}-2^\frac{k}{2}}$, else ${\Delta g_k = 2^{k+1}-2^\frac{k+1}{2}}$ so the second equality holds trivially. Moreover, since ${p\geq 3}$, the first one follows from Bounds (\ref{nbratplgsr}) and (\ref{nbpldeg12gsr}) which gives ${N_k \geq 2^{k+2}}$.
\qed
\end{Proof}

\begin{lemm}\label{lemme_bornesupgsr}
Let $H_k$ be a step of one of the towers $T/\F_{p^2}$ or $T/\F_p$. One has:
$$
\sup \big \{ n \in \N \; | \;  N_k \geq 2n +2g_k -1\big \} \geq 2^{k}(p-3)+2.
$$
\end{lemm}

\begin{Proof}
From Bounds (\ref{nbratplgsr}) and (\ref{nbpldeg12gsr}) for $N_k$ and Lemma \ref{lemme_genregsr} i), we get 
\begin{eqnarray*}
N_k - 2g_k +1 & \geq & 2^{k+1}(p-1) -2(2^{k+1}-2\cdot2^\frac{k+1}{2}+1) +1\\
		     & = & 2^{k+1}(p-3) + 4\cdot2^\frac{k+1}{2} - 1 \\
		     & \geq & 2^{k+1}(p-3) + 4 \mbox{ since } k\geq0.
\end{eqnarray*}\qed
\end{Proof}

\subsection{General results for $\mu_q(n)$}
In \cite{balb}, Ballet and Le Brigand proved the following useful result:
\begin{theo}\label{existdivnonspe}
Let $F/\F_q$ be an algebraic function field of genus $g\geq 2$. If $q\geq4$, then there exists a non-special divisor of degree $g-1$.
\end{theo}

The four following lemmas prove the existence of a "good" step of the towers defined in Section \ref{sectdeftowers}, that is to say a step that will be optimal for the bilinear complexity of multiplication:

\begin{lemm}\label{lemme_placedegn}
Let $n \geq \frac{1}{2}\left(q^2+1+\epsilon(q^2)\right)$ be an integer. If $q=p^r\geq4$, then there exists a step $F_{k,s}/\F_{q^2}$ of the tower $T_2/\F_{q^2}$ such that all the three following conditions are verified:
\begin{enumerate}[(1)]
	\item there exists a non-special divisor of degree $g_{k,s}-1$ in $F_{k,s}/\F_{q^2}$,
	\item there exists a place of $F_{k,s}/\F_{q^2}$ of degree $n$,
	\item $N_1(F_{k,s}/\F_{q^2}) \geq 2n + 2g_{k,s}-1$.
\end{enumerate}
Moreover, the first step for which both Conditions (2) and (3) are verified is the first step for which (3) is verified.
\end{lemm}

\begin{Proof}
Note that $n \geq 9$ since $q\geq4$ and ${n \geq \frac{1}{2}(q^2+1) \geq 8.5}$. Fix $1 \leq k \leq n-4$ and ${s \in \{0, \ldots, r\}}$. First, we prove that Condition (2) is verified. Lemma~\ref{lemme_genre}~iv) gives:
\begin{eqnarray}
	\nonumber 2g_{k,s}+1 & \leq & 2\frac{q^k(q+1)-q^\frac{k}{2}(q-1)}{p^{r-s}} +1\\
	\nonumber		& = & 2p^s\left(q^{k-1}(q+1)-q^\frac{k}{2}\frac{q-1}{q}\right) +1\\
					& \leq & 2q^{k-1}p^s(q+1) \ \  \ \mbox{\ since } 2p^sq^\frac{k}{2}\frac{q-1}{q}\geq1 \label{eqgenre1}\\
	\nonumber		& \leq & 2q^k(q^2-1).
\end{eqnarray}
On the other hand, one has ${n-1 \geq k+3 > k+\frac{1}{2}+2}$ so $n-1 \geq \log_q(q^k)+\log_q(2)+\log_q(q+1)$. This gives ${q^{n-1} \geq 2q^k(q+1)}$, hence $q^{n-1}(q-1) \geq 2q^k(q^2-1)$. Therefore, one has ${2g_{k,s}+1 \leq q^{n-1}(q-1)}$ which ensure us that Condition (2) is satisfied according to Corollary 5.2.10 in \cite{stic}.\\
Now suppose also that ${k \geq \log_q\left(\frac{2n}{5}\right)+1}$. Note that for all $n\geq 9$ there exists such an integer $k$ since the size of the interval $[\log_q\left(\frac{2n}{5}\right)+1 , n-4]$ is bigger than ${9-4-\log_4\left(\frac{2\cdot9}{5}\right)-1 \geq 3 >1}$. Moreover such an integer $k$ verifies ${q^{k-1} \geq \frac{2}{5}n}$, so ${n \leq \frac{1}{2}q^{k-1}(q+1)(q-3)}$ since $q\geq4$. Then one has
\begin {eqnarray*}
	2n+2g_{k,s}-1 & \leq & 2n+2g_{k,s}+1\\
			& \leq & 2n + 2q^{k-1}p^s(q+1) \ \  \ \mbox{\ according to (\ref{eqgenre1})}\\
			& \leq & q^{k-1}(q+1)(q-3) + 2q^{k-1}p^s(q+1)\\
			& \leq  & q^{k-1}p^s(q+1)(q-1) \\
			& = & (q^2-1)q^{k-1}p^s
\end{eqnarray*}
which gives ${N_1(F_{k,s}/\F_{q^2}) \geq 2n + 2g_{k,s}-1}$ according to Proposition \ref{subfield} (3). Hence, for any integer $k \in [\log_q\left(\frac{2n}{5}\right)+1 , n-4]$, Conditions (2) and (3) are satisfied and the smallest integer $k$ for which they are both satisfied is the smallest integer $k$ for which Condition (3) is satisfied.\\
To conclude, remark that for such an integer $k$, Condition (1) is easily verified from Theorem \ref{existdivnonspe} since $q\geq 4$ and ${g_{k,s} \geq g_2\geq 6}$ according to Formula (\ref{genregs}).\\
\qed

\end{Proof}

This is a similar result for the tower $T_3/\F_q$:

\begin{lemm}\label{lemme_placedegn2}
Let $n \geq \frac{1}{2}\left(q+1+\epsilon(q)\right)$ be an integer. If $q=p^r\geq 4$, then there exists a step $G_{k,s}/\F_q$ of the tower $T_3/\F_q$ such that all the three following conditions are verified:
\begin{enumerate}[(1)]
	\item there exists a non-special divisor of degree $g_{k,s}-1$ in $G_{k,s}/\F_q$,
	\item there exists a place of $G_{k,s}/\F_q$ of degree $n$, 
	\item $N_1(G_{k,s}/\F_q)+2N_2(G_{k,s}/\F_q) \geq 2n + 2g_{k,s}-1$.
\end{enumerate}
Moreover, the first step for which both Conditions (2) and (3) are verified is the first step for which (3) is verified.
\end{lemm}

\begin{Proof} Note that $n \geq 5$ since $q\geq4$, ${\epsilon(q) \geq \epsilon(4)=4}$ and ${n \geq \frac{1}{2}(q+1+\epsilon(q)) \geq 4.5}$.
First, we focus on the case $n\geq13$. Fix $1 \leq k \leq \frac{n-7}{2}$ and ${s \in \{0, \ldots, r\}}$. One has ${2p^sq^k\frac{q+1}{\frac{\sqrt{q}}{2}} \leq q^\frac{n-1}{2}}$ since
$${
\frac{n-1}{2} \geq k + 3 = k -\frac{1}{2} +1+1+\frac{3}{2} \geq \log_q(q^{k-\frac{1}{2}}) +\log_q(4)+\log_q(p^s)+\log_q(q+1)}.
$$
Hence ${2p^sq^k(q+1) \leq q^\frac{n-1}{2}(\sqrt{q}-1)}$ since ${\frac{\sqrt{q}}{2}\leq \sqrt{q}-1}$ for $q\geq4$. According to (\ref{eqgenre1}) in the previous proof, this proves that Condition (2) is satisfied.\\
The same reasoning as in the previous proof shows that Condition (3) is also satisfied as soon as ${k \geq \log_q\left(\frac{2n}{5}\right)+1}$. Moreover, for $n\geq13$, the interval $[\log_q\left(\frac{2n}{5}\right)+1 , \frac{n-7}{2}]$ contains at least one integer and the smallest integer $k$ in this interval is the smallest integer $k$ for which Condition (3) is verified. Furthermore, for such an integer $k$, Condition (1) is easily verified from Theorem \ref{existdivnonspe} since $q\geq 4$ and ${g_{k,s} \geq g_2\geq 6}$ according to Formula (\ref{genregs}).

To complete the proof, we want to focus on the case $5\leq n\leq12$. For this case, we have to look at the values of $q=p^r$ and $n$ for which we have both ${n\geq \frac{1}{2}\left(q+1+\epsilon(q)\right)}$ and ${5 \leq n \leq 12}$. For each value of $n$ such that these two inequalities are satisfied, we have to check that Conditions (1), (2) and (3) are verified. In this aim, we use the KASH packages \cite{kash} to compute the genus and number of places of degree 1 and 2 of the first steps of the tower $T_3/\F_q$. Thus we determine the first step $G_{k,s}/\F_q$ that satisfied all the three Conditions (1), (2) and (3). We resume our results in the following table:
\hspace{-8em} 
$$\begin{array}{|c|c|c|c|}
	\hline
	q=p^r & 2^2 & 2^3 & 3^2 \\
	 \hline
 	\epsilon(q) & 4 & 5 & 6 \\
	\hline
	\frac{1}{2}\left(q+1+\epsilon(q)\right) & 4.5 & 7 & 8  \\
	\hline
	n \hbox{ to be considered} & 5 \leq n \leq 12 & 7 \leq n \leq 12 & 8 \leq n \leq 12 \\
	\hline
	(k,s) & (1,1) & (1,1) & (1,1) \\
 	\hline
 	N_1(G_{k,s}/\F_q) & 5 & 9 & 10  \\
	\hline
	N_2(G_{k,s}/\F_q) & 14 & 124 & 117 \\
	\hline
	\Gamma(G_{k,s}/\F_q) & 15 & 117 & 113 \\
	\hline
	g_{k,s} & 2 & 12 & 9 \\
	\hline
	2g_{k,s}+1 & 5 & 25 & 19 \\
	\hline
	q^\frac{n-1}{2}(\sqrt{q}-1) \geq \ldots & 16 & 936 & 4374 \\
	\hline
\end{array}
$$

\vspace{-0.5em}
\hspace{-8em} 
$$
\begin{array}{|c|c|c|c|c|}
	\hline
	q=p^r & 5 & 7 & 11 & 13 \\
	 \hline
 	\epsilon(q) & 4 & 5 & 6 & 7\\
	\hline
	\frac{1}{2}\left(q+1+\epsilon(q)\right) & 5 & 6.5 & 9 & 10.5 \\
	\hline
	n \hbox{ to be considered} &  5 \leq n \leq 12 & 7 \leq n \leq 12 & 9 \leq n \leq 12 & 11 \leq n \leq 12 \\
	\hline
	(k,s) & (2,0) & (2,0) & (2,0) & (2,0) \\
 	\hline
 	N_1(G_{k,s}/\F_q) & 6 & 8 & 12 & 14 \\
	\hline
	N_2(G_{k,s}/\F_q) & 60 & 168 & 660 & 1092 \\
	\hline
	\Gamma(G_{k,s}/\F_q) & 53 & 151.5  & 611.5  & 1021.5  \\
	\hline
	g_{k,s} & 10 & 21 & 55 & 78 \\
	\hline
	2g_{k,s}+1 & 21 & 43 & 11 & 157 \\
	\hline
	q^\frac{n-1}{2}(\sqrt{q}-1) \geq \ldots & 30 & 564 & 33917 & 967422 \\
	\hline
\end{array}
$$
\vspace{1em}

In this table, one can check that for each value of $q$ and $n$ to be considered and every corresponding step 
$G_{k,s}/\F_q$ one has simultaneously:
\begin{itemize}
	\item $g_{k,s}\geq2$ so Condition (1) is verified according to Theorem \ref{existdivnonspe},
	\item $2g_{k,s}+1 \leq q^\frac{n-1}{2}(\sqrt{q}-1)$ so Condition (2) is verified.
	\item $\Gamma(G_{k,s}/\F_q):=\frac{1}{2}\left(N_1(G_{k,s}/\F_q)+2N_2(G_{k,s}/\F_q)-2g_{k,s}+1\right) \geq n$ so Condition (3) is verified.
\end{itemize}

\qed
\end{Proof}

This is a similar result for the tower $T/\F_{p^2}$:

\begin{lemm}\label{lemme_placedegngsr}
Let $p\geq5$ and $n \geq \frac{1}{2}\left(p^2+1+\epsilon(p^2)\right)$. There exists a step $H_k/\F_{p^2}$ of the tower $T/\F_{p^2}$ such that the three following conditions are verified:
\begin{enumerate}[(1)]
	\item there exists a non-special divisor of degree $g_k-1$ in $H_k/\F_{p^2}$,
	\item there exists a place of $H_k/\F_{p^2}$ of degree $n$,
	\item $N_1(H_k/\F_{p^2}) \geq 2n + 2g_k - 1$.
\end{enumerate}
Moreover the first step for which all the three conditions are verified is the first step for which (3) is verified.
\end{lemm}

\begin{Proof} Note that ${n \geq \frac{1}{2}(5^2+1+\epsilon(5^2)) = 18}$. We first prove that for all integers $k$ such that ${2 \leq k \leq n - 2}$, we have  ${2g_k+1 \leq p^{n-1}(p-1)}$ , so Condition (2) is verified according to Corollary 5.2.10 in \cite{stic2}. Indeed, for such an integer $k$, since ${p\geq5}$ one has ${k \leq \log_2(p^{n-2}) \leq  \log_2(p^{n-1}-1)}$, thus $k+2 \leq \log_2\left(4(p^{n-1}-1)\right) \leq \log_2 (4p^{n-1}-1)$ and it follows that ${2^{k+2}+1 \leq 4p^{n-1}}$. Hence ${2\cdot2^{k+1} +1 \leq p^{n-1}(p-1)}$ since ${p\geq5}$, which gives the result according to Lemma \ref{lemme_genregsr} ii).\\
We prove now that for ${k\geq \log_2 (2n-1)-2}$, Condition (3) is verified. Indeed, for such an integer $k$, we have ${k +2\geq \log_2 (2n-1)}$, so ${2^{k+2} \geq 2n-1}$. Hence we get ${2^{k+3} \geq 2n+2^{k+2}-1}$ and so ${2^{k+1}(p-1)\geq 2^{k+1}\cdot4 \geq 2n+2^{k+2}-1}$ since ${p\geq5}$. Thus we have ${N_1(H_k/\F_{p^2}) \geq 2n + 2g_k - 1}$ according to Bound (\ref{nbratplgsr}) and Lemma~\ref{lemme_genregsr}~ii).\\
Hence, we have proved that for any integers ${n\geq 18}$ and ${k \geq 2}$ such that\linebreak[4] ${\log_2 (2n-1)-2 \leq k \leq  n - 2}$, both Conditions (2) and (3) are verified. Moreover, note that for any ${n\geq 18}$, there exists an integer $k \geq 2$ in the interval ${\big[ \log_2 (2n-1)-2; n - 2 \big]}$. Indeed, ${\log_2 (2\cdot 18-1)-2 \simeq 3.12 >2}$ and the size of this interval increases with $n$ and is greater than 1 for $n=18$. To conclude, remark that for such an integer $k$, Condition (1) is easily verified from Theorem \ref{existdivnonspe} since $p^2\geq 4$ and ${g_k \geq g_2=3}$ according to Formula (\ref{genregsr}).\\
\qed
\end{Proof}

This is a similar result for the tower $T/\F_p$:

\begin{lemm}\label{lemme_placedegngsr2}
Let $p\geq5$ and $n \geq \frac{1}{2}\left(p+1+\epsilon(p)\right)$. There exists a step $H_k/\F_p$ of the tower $T/\F_p$ such that  the three following conditions are verified:
\begin{enumerate}[(1)]
	\item there exists a non-special divisor of degree $g_k-1$ in $H_k/\F_p$,
	\item there exists a place of $H_k/\F_p$ of degree $n$,
	\item $N_1(H_k/\F_p) + 2N_1(H_k/\F_p) \geq 2n + 2g_k - 1$.
\end{enumerate}
Moreover the first step for which all the three conditions are verified is the first step for which (3) is verified.
\end{lemm}

\begin{Proof} Note that ${n \geq \frac{1}{2}(5+1+\epsilon(5)) =5}$.
We first prove that for all integers $k$ such that ${2 \leq k \leq n - 3}$, we have  ${2g_k+1 \leq p^\frac{n-1}{2}(\sqrt{p}-1)}$,
 so Condition (2) is verified according to Corollary 5.2.10 in \cite{stic2}.
Indeed, for such an integer $k$, since ${p\geq5}$ and ${n\geq5}$ 
one has ${\log_2(p^\frac{n-1}{2}-1) \geq \log_2(5^\frac{n-1}{2}-1) \geq \log_2(2^{n-1}) = n-1}$.
Thus ${k+2\leq n-1\leq\log_2(p^\frac{n-1}{2}-1)}$ and it follows from Lemma \ref{lemme_genregsr} ii) that \linebreak[4]${2g_k+1 \leq 2^{k+2}+1 \leq p^\frac{n-1}{2} \leq p^\frac{n-1}{2}(\sqrt{p}-1)}$, 
which gives the result.\\
The same reasoning as in the previous proof shows that Condition (3) is also satisfied as soon as  ${k\geq \log_2 (2n-1)-2}$.
Hence, we have proved that for any integers ${n\geq 5}$ and ${k \geq 2}$ such that ${\log_2 (2n-1)-2 \leq k \leq  n - 3}$, 
both Conditions (2) and (3) are verified.
Moreover, note that  the size of the interval ${\big[ \log_2 (2n-1)-2; n - 3 \big]}$
increases with $n$ and that for any ${n\geq 5}$, this interval contains at least one integer $k \geq 2$.
To conclude, remark that for such an integer $k$, Condition (1) is easily verified from Theorem \ref{existdivnonspe} since $p\geq 4$ and ${g_k \geq g_2=3}$ according to Formula (\ref{genregsr}).\\
\qed
\end{Proof}

Now we establish general bounds for the bilinear complexity of multiplication by using derivative evaluations on places of degree one (respectively places of degree one and two).

\begin{theo}{\label{thm_arnaud1}}
Let $q$ be a prime power and $n>1$ be an integer. If there exists an algebraic function field $F/ \F_q$ of genus $g$ with $N$ places of degree 1 and an integer $0 < a \leq N$ such that 
\begin{enumerate}[(i)]
	\item there exists $\mathcal{R}$, a non-special divisor of degree $g-1$,
	\item there exists $Q$, a place of degree $n$,
	\item $N+a \geq 2n+2g-1$.
\end{enumerate}
Then
$$
\mu_q(n) \leq 2n +g-1+a\mbox{.}
$$
\end{theo}

\begin{Proof} Let $\mathcal{P}:=\{P_1, \ldots, P_N\}$ be a set of $N$ places of degree 1 and $\mathcal{P}'$ be a subset of $\mathcal{P}$ with cardinal number $a$. According to Lemma 2.7 in \cite{bapi}, we can choose an effectif divisor $\mathcal{D}$ equivalent to $Q+\mathcal{R}$ such that ${\mathrm{supp}(\mathcal{D}) \cap \mathcal{P} = \varnothing}$. We define the maps $Ev_Q$ and $Ev_\mathcal{P}$ as in Theorem \ref{theo_evalder} with $u_i=2$ if $P_i \in \mathcal{P}'$ and $u_i=1$ if $P_i \in \mathcal{P}\backslash\mathcal{P}'$. Then $Ev_Q$ is bijective, since $\ker Ev_Q = \mathcal{L}(\mathcal{D}-Q)$ with ${\dim(\mathcal{D}-Q) = \dim(R) =0}$ and ${\dim (\mathrm{im}\, Ev_Q )= \dim \mathcal{D} = \deg\mathcal{D} -g +1 + \mathrm{i}(\mathcal{D}) \geq n}$ according to Riemann-Roch Theorem. Thus $\dim (\mathrm{im}\, Ev_Q ) =n$. Moreover, $Ev_\mathcal{P}$ is injective. Indeed, ${\ker Ev_\mathcal{P} = \mathcal{L}(2\mathcal{D}-\sum_{i=1}^N u_iP_i)}$ with $\deg (2\mathcal{D}-\sum_{i=1}^N u_iP_i) =2(n+g-1)-N-a <0$. Furthermore, one has $\mathrm{rk}\, Ev_\mathcal{P} = \dim(2\mathcal{D})= \deg(2\mathcal{D})-g+1+\mathrm{i}(2\mathcal{D})$, and $\mathrm{i}(2\mathcal{D})=0$ since $2\mathcal{D} \geq \mathcal{D} \geq \mathcal{R}$ with ${\mathrm{i}(\mathcal{R})=0}$. So ${\mathrm{rk}\, Ev_\mathcal{P} = 2n+g-1}$, and we can extract a subset $\mathcal{P}_1$ from $\mathcal{P}$ and a subset $\mathcal{P}_1'$ from $\mathcal{P}'$ with cardinal number $N_1\leq N$ and $a_1\leq a$, such that:
\begin{itemize}
	\item $N_1+a_1 = 2n+g-1$,
	\item the map $Ev_{\mathcal{P}_1}$ defined as $Ev_\mathcal{P}$ with $u_i=2$ if $P_i \in \mathcal{P}_1'$ and $u_i=1$ if $P_i \in \mathcal{P}_1\backslash\mathcal{P}_1'$, is injective.
\end{itemize}
According to Theorem \ref{theo_evalder}, this leads to $\mu_q(n) \leq N_1+2a_1 \leq N_1+a_1+a$ which gives the result.
\qed
\end{Proof}

\begin{theo}{\label{thm_arnaud2}}
Let $q$ be a prime power and $n>1$ be an integer. If there exists an algebraic function field $F/ \F_q$ of genus $g$ with $N_1$ places of degree 1, $N_2$ places of degree 2 and two integers  $0 < a_1 \leq N_1$, $0 < a_2 \leq N_2$ such that 
\begin{enumerate}[(i)]
	\item there exists  $\mathcal{R}$, a non-special divisor of degree $g-1$,
	\item there exists $Q$, a place of degree $n$,
	\item $N_1+a_1 +2(N_2+a_2) \geq 2n+2g-1$.
\end{enumerate}
Then
$$
\mu_q(n) \leq 2n + g +N_2 + a_1 + 4a_2
$$
and
$$
\mu_q(n) \leq 3n+\frac{3}{2}g+\frac{a_1}{2}+3a_2. 
$$
\end{theo}

\begin{Proof}
Let $\mathcal{P}_1:=\{P_1, \ldots, P_{N_1}\}$ be a set of $N_1$ places of degree 1 and $\mathcal{P}_1'$ be a subset of $\mathcal{P}_1$ with cardinal number $a_1$. Let $\mathcal{P}_2:=\{Q_1, \ldots, Q_{N_2}\}$ be a set of $N_2$ places of degree 2 and $\mathcal{P}_2'$ be a subset of $\mathcal{P}_2$ with cardinal number~$a_2$. According to Lemma 2.7 in \cite{bapi}, we can choose an effectif divisor $\mathcal{D}$ equivalent to $Q+\mathcal{R}$ such that ${\mathrm{supp}(\mathcal{D}) \cap (\mathcal{P}_1\cup \mathcal{P}_2) = \varnothing}$. We define the maps $Ev_Q$ and $Ev_\mathcal{P}$ as in Theorem \ref{theo_evalder} with $u_i=2$ if $P_i \in \mathcal{P}_1' \cup \mathcal{P}_2'$ and $u_i=1$ if ${P_i \in (\mathcal{P}_1\backslash\mathcal{P}_1') \cup (\mathcal{P}_2\backslash\mathcal{P}_2')}$. Then the same raisoning as in the previous proof shows that $Ev_Q$ is bijective. Moreover, $Ev_\mathcal{P}$ is injective. Indeed, ${\ker Ev_\mathcal{P} = \mathcal{L}(2\mathcal{D}-\sum_{i=1}^N u_iP_i)}$ with $\deg (2\mathcal{D}-\sum_{i=1}^N u_iP_i) =2(n+g-1)-(N_1+a_1+2(N_2+a_2)) <0$. Furthermore, one has $\mathrm{rk}\, Ev_\mathcal{P} = \dim(2\mathcal{D})= \deg(2\mathcal{D})-g+1+\mathrm{i}(2\mathcal{D})$, and $\mathrm{i}(2\mathcal{D})=0$ since $2\mathcal{D} \geq \mathcal{D} \geq \mathcal{R}$ with ${\mathrm{i}(\mathcal{R})=0}$. So ${\mathrm{rk}\, Ev_\mathcal{P} = 2n+g-1}$, and we can extract a subset $\Tilde{\mathcal{P}}_1$ from $\mathcal{P}_1$, a subset $\Tilde{\mathcal{P}}_1'$ from $\mathcal{P}_1'$, a subset $\Tilde{\mathcal{P}}_2$ from $\mathcal{P}_2$ and a subset $\Tilde{\mathcal{P}}_2'$ from $\mathcal{P}_2'$  with respective cardinal numbers $\Tilde{N}_1\leq N_1$,  $\Tilde{a}_1\leq a_1$, $\Tilde{N}_2\leq N_2$ and  $\Tilde{a}_2\leq a_2$, such that:
\begin{itemize}
	\item $2n + g \geq \Tilde{N}_1+\Tilde{a}_1 +2(\Tilde{N}_2+\Tilde{a}_2) \geq 2n+g-1$,
	\item the map $Ev_{\Tilde{\mathcal{P}}}$ defined as $Ev_\mathcal{P}$ with $u_i=2$ if $P_i \in \Tilde{\mathcal{P}}_1'\cup \Tilde{\mathcal{P}}_2'$ and $u_i=1$ if $(\Tilde{\mathcal{P}}_1\backslash \Tilde{\mathcal{P}}_1') \cup (\Tilde{\mathcal{P}}_2\backslash \Tilde{\mathcal{P}}_2')$, is injective.
\end{itemize}
According to Theorem \ref{theo_evalder}, this leads to $\mu_q(n) \leq \Tilde{N}_1+2\Tilde{a}_1 + 3(\Tilde{N}_2 +2\Tilde{a}_2)$ since $M_k(2) \leq3$ for all prime power $k$. Hence, one has the first result since \linebreak[4]${\Tilde{N}_1+\Tilde{a}_1 + 2(\Tilde{N}_2 +\Tilde{a}_2)\leq 2n+g}$ and the second one since ${\frac{\Tilde{a}_1}{2}+\Tilde{N}_2+\Tilde{a}_2 \leq \frac{g}{2}+n}$.
\qed
\end{Proof}

\subsection{New upper bounds for $\mu_{q}(n)$} \label{sectbornesarnaud}

Here, we give a detailed proof of Bound (i) of Theorem \ref{bornes_arnaud1} and we give an improvement of Bound (ii).
Moreover, we correct the bound for $\mu_{p^2}(n)$ given in \cite{arna1} and ameliorate the unproved bound for $\mu_p(n)$. Namely, we prove:

\begin{theo}\label{theo_arnaud1}
Let ${q=p^r \geq 4}$ be a power of the prime $p$. Then
	\begin{equation*}
	\begin{array}{l}
	\mbox{(i) \ If ${q=p^r \geq 4}$, then  }\mu_{q^2}(n) \leq 2 \left(1 + \frac{p}{q-3 + (p-1)\left(1- \frac{1}{q+1}\right)} \right)n,\\
	\\
	\mbox{(ii) \ If ${q=p^r \geq 4}$, then }\mu_{q}(n) \leq 3 \left(1 + \frac{p}{q-3 + (p-1)\left(1- \frac{1}{q+1}\right)} \right)n.\\
	\\
	\mbox{(iii) \ If $p\geq5$, then  }\mu_{p^2}(n) \leq 2 \left(1 + \frac{2}{p-\frac{33}{16}} \right)n.\\
	\\
	\mbox{(iv) \ If $p\geq5$, then  }\mu_{p}(n) \leq 3\left(1 + \frac{2}{p-\frac{33}{16}} \right)n.
	\end{array}
	\end{equation*}

\end{theo}

\begin{Proof}
\begin{enumerate}[(i)]
	\item Let $n\geq \frac{1}{2}(q^2+1+{\epsilon (q^2) })$. Otherwise, we already know from Theorems \ref{thm_wdg} and  \ref{thm_shokr} that $\mu_{q^2}(n) \leq 2n$. According to Lemma \ref{lemme_placedegn}, there exists a step of the tower $T_2/\F_{q^2}$ on which we can apply Theorem \ref{thm_arnaud1} with $a=0$. We denote by $F_{k,s+1}/\F_{q^2}$ the first step of the tower that suits the hypothesis of Theorem \ref{thm_arnaud1} with $a=0$, i.e. $k$ and $s$ are integers such that ${N_{k,s+1} \geq 2n+2g_{k,s+1}-1}$ and ${N_{k,s}< 2n+2g_{k,s}-1}$, where ${N_{k,s}:=N_1(F_{k,s}/\F_{q^2})}$ and ${g_k:=g(F_{k,s})}$. We denote by $n_0^{k,s}$ the biggest integer such that ${N_{k,s}\geq 2n_0^{k,s}+2g_{k,s}-1}$, i.e. ${n_0^{k,s} = \sup \big\{n \in \N \, \vert \, 2n \leq N_{k,s}-2g_{k,s}+1\big\}}$. To perform multiplication in $\F_{q^{2n}}$, we have the following alternative:
	\begin{enumerate}[(a)]
		\item use the algorithm on the step $F_{k,s+1}$. In this case, a bound for the bilinear complexity is given by Theorem \ref{thm_arnaud1} applied with $a=0$:
		$$
		\mu_{q^2}(n) \leq 2n+g_{k,s+1}-1= 2n+g_{k,s}-1 +\Delta g_{k,s}.
		$$
		(Recall that $\Delta g_{k,s} := g_{k,s+1} - g_{k,s}$)
		\item use the algorithm on the step $F_{k,s}$ with an appropriate number of derivative evaluations. Let $a:= 2(n-n_0^{k,s})$ and suppose that $a \leq N_{k,s}$. Then ${N_{k,s} \geq 2n_0^{k,s}+2g_{k,s}-1}$ implies that ${N_{k,s} +a \geq 2n+2g_{k,s}-1}$ so Condition (iii) of Theorem \ref{thm_arnaud1} is satisfied. Thus, we can perform $a$ derivative evaluations in the algorithm using the step $F_{k,s}$ and we have:
		$$
		\mu_{q^2}(n) \leq 2n+g_{k,s}-1+a.
		$$
	\end{enumerate}
	Thus, if $a \leq N_{k,s}$ Case (b) gives a better bound as soon as ${a<\Delta g_{k,s}}$. Since we have from Lemma \ref{lemme_delta} both ${N_{k,s} \geq D_{k,s}}$ and ${\Delta g_{k,s} \geq D_{k,s}}$, if ${a\leq D_{k,s}}$ then we can perform $a$ derivative evaluations on places of degree 1 in the step $F_{k,s}$ and Case (b) gives a better bound then Case (a).\\
	For $x \in \mathbb{R}^{+}$ such that ${N_{k,s+1} \geq 2[x]+2g_{k,s+1}-1}$ and ${N_{k,s} < 2[x]+2g_{k,s}-1}$, we define the function $\Phi_{k,s}(x)$ as follow:
	$$
	\Phi_{k,s}(x) = \left\{\begin{array}{ll}
				2x+g_{k,s}-1+2(x-n_0^{k,s}) & \mbox{if } 2(x- n_0^{k,s}) < D_{k,s}\\
				2x+g_{k,s+1}-1& \mbox{else}.
			      \end{array} \right.
	$$
	We define the function $\Phi$ for all ${x\geq0}$ as the minimum of the functions $\Phi_{k,s}$ for which $x$ is in the domain of $\Phi_{k,s}$. This function is piecewise linear with two kinds of piece: those which have slope $2$ and those which have slope~$4$. Moreover, since the y-intercept of each piece grows with $k$ and $s$, the graph of the function $\Phi$ lies below any straight line that lies above all the points ${\big(n_0^{k,s}+\frac{D_{k,s}}{2}, \Phi(n_0^{k,s}+\frac{D_{k,s}}{2})\big)}$, since these are the \textit{vertices} of the graph. Let ${X:=n_0^{k,s}+\frac{D_{k,s}}{2}}$, then
	\begin{eqnarray*}
	\Phi(X) & \leq & 2X + g_{k,s+1} -1\\
		    & \leq & 2X+ g_{k,s+1}\\
		    & = & 2\left(1 + \frac{g_{k,s+1}}{2X}\right)X.
	\end{eqnarray*}
We want to give a bound for $\Phi(X)$ which is independent of $k$ and $s$.

Recall that $D_{k,s} :=(p-1)p^sq^k$, and
$$
2n_0^{k,s}  \geq q^{k-1}p^s(q+1)(q-3) \  \  \  \mbox{by Lemma \ref{lemme_bornesup}}
$$
and
$$
g_{k,s+1} \leq q^{k-1}(q+1)p^{s+1}  \  \  \  \mbox{by Lemma \ref{lemme_genre} (iii).}
$$
So we have
\begin{eqnarray*}
	\frac{g_{k,s+1}}{2X} & = & \frac{g_{k,s+1}}{2n_0^{k,s}+D_{k,s}} \\
				       & \leq & \frac{q^{k-1}(q+1)p^{s+1}}{q^{k-1}p^s(q+1)(q-3) + (p-1)p^sq^k } \\
		    		       & = & \frac{q^{k-1}(q+1)p^sp}{q^{k-1}(q+1)p^s\left(q-3 + (p-1)\frac{q}{q+1}\right) }\\
		    		       & = & \frac{p}{(q-3)+(p-1)\frac{q}{q+1}}
\end{eqnarray*}
Thus, the graph of the function $\Phi$ lies below the line ${y=2\left(1 + \frac{p}{(q-3)+(p-1)\frac{q}{q+1}}\right)x}$. In particular, we get
$$
\Phi(n) \leq 2\left(1 + \frac{p}{(q-3)+(p-1)\frac{q}{q+1}}\right)n.
$$
	\item Let $n\geq \frac{1}{2}(q+1+{\epsilon (q) })$. Otherwise, we already know from Theorems \ref{thm_wdg} and  \ref{thm_shokr} that $\mu_{q}(n) \leq 2n$. According to Lemma \ref{lemme_placedegn2}, there exists a step of the tower $T_3/\F_{q}$ on which we can apply Theorem \ref{thm_arnaud2} with $a_1=a_2=0$. We denote by $G_{k,s+1}/\F_{q}$ the first step of the tower that suits the hypothesis of Theorem \ref{thm_arnaud2} with $a_1=a_2=0$, i.e. $k$ and $s$ are integers such that ${N_{k,s+1} \geq 2n+2g_{k,s+1}-1}$ and ${N_{k,s}< 2n+2g_{k,s}-1}$, where \linebreak[4]${N_{k,s}:=N_1(G_{k,s}/\F_q)+2N_2(G_{k,s}/\F_q)}$ and ${g_{k,s}:=g(G_{k,s})}$. We denote by $n_0^{k,s}$ the biggest integer such that ${N_{k,s}\geq 2n_0^{k,s}+2g_{k,s}-1}$, i.e. \linebreak[4] ${n_0^{k,s} = \sup \big\{n \in \N \, \vert \, 2n \leq N_{k,s}-2g_{k,s}+1\big\}}$. To perform multiplication in $\F_{q^n}$, we have the following alternative:
	\begin{enumerate}[(a)]
		\item use the algorithm on the step $G_{k,s+1}$. In this case, a bound for the bilinear complexity is given by Theorem \ref{thm_arnaud2} applied with $a_1=a_2=0$:
		$$
		\mu_q(n) \leq 3n+\frac{3}{2}g_{k,s+1}= 3n_0^{k,s}+\frac{3}{2}g_{k,s}+3(n-n_0^{k,s}) + \frac{3}{2}\Delta g_{k,s}.
		$$
		\item use the algorithm on the step $G_{k,s}$ with an appropriate number of derivative evaluations. Let ${a_1+2a_2:= 2(n-n_0^{k,s})}$ and suppose that ${a_1+2a_2 \leq N_{k,s}}$. Then ${N_{k,s}\geq 2n_0^{k,s}+2g_{k,s}-1}$ implies that ${N_{k,s} +a_1+2a_2 \geq 2n+2g_{k,s}-1}$. Thus we can perform $a_1+a_2$ derivative evaluations in the algorithm using the step $G_{k,s}$ and we have:
		$$
		\mu_q(n) \leq 3n+\frac{3}{2}g_{k,s}+\frac{3}{2}(a_1+2a_2)=3n_0^{k,s}+\frac{3}{2}g_{k,s}+6(n-n_0^{k,s}).
		$$
	\end{enumerate}
	Thus, if $a_1+2a_2 \leq N_{k,s}$ Case (b) gives a better bound as soon as ${n-n_0^{k,s}<\frac{1}{2}\Delta g_{k,s}}$. Since we have from Lemma \ref{lemme_delta} both $N_{k,s} \geq D_{k,s}$ and $\frac{1}{2}\Delta g_{k,s} \geq \frac{1}{2}D_{k,s}$, if $a_1+2a_2\leq D_{k,s}$, i.e. $n-n_0^{k,s} \leq \frac{1}{2} D_{k,s}$, then we can perform $a_1$ derivative evaluations on places of degree 1 and $a_2$ derivative evaluations on places of degree 2 in the step $G_{k,s}$ and Case (b) gives a better bound then Case (a).\\
	For $x \in \mathbb{R}^{+}$ such that ${N_{k,s+1} \geq 2[x]+2g_{k,s+1}-1}$ and ${N_{k,s} < 2[x]+2g_{k,s}-1}$, we define the function $\Phi_{k,s}(x)$ as follow:
	$$
	\Phi_{k,s}(x) = \left\{\begin{array}{ll}
				3x+\frac{3}{2}g_{k,s}+3(x-n_0^{k,s}) & \mbox{if } x- n_0^{k,s} < \frac{D_{k,s}}{2}\\
				& \\
				3x+\frac{3}{2}g_{k,s+1}& \mbox{else}.
			      \end{array} \right.
	$$
	We define the function $\Phi$ for all ${x\geq0}$ as the minimum of the functions $\Phi_{k,s}$ for which $x$ is in the domain of $\Phi_{k,s}$. This function is piecewise linear with two kinds of piece: those which have slope $3$ and those which have slope~$6$. Moreover, since the y-intercept of each piece grows with $k$ and $s$, the graph of the function $\Phi$ lies below any straight line that lies above all the points ${\big(n_0^{k,s}+\frac{D_{k,s}}{2}, \Phi(n_0^{k,s}+\frac{D_{k,s}}{2})\big)}$, since these are the \textit{vertices} of the graph. Let ${X:=n_0^{k,s}+\frac{D_{k,s}}{2}}$, then
	\begin{eqnarray*}
	\Phi(X) & \leq & 3X + \frac{3}{2}g_{k,s+1} \\
		    & = & 3\left(1 + \frac{g_{k,s+1}}{2X}\right)X.
	\end{eqnarray*}
We want to give a bound for $\Phi(X)$ which is independent of $k$ and $s$.

Recall that $D_{k,s} :=(p-1)p^sq^k$, and
$$
n_0^{k,s}  \geq \frac{1}{2}q^{k-1}p^s(q+1)(q-3) \  \  \  \mbox{by Lemma \ref{lemme_bornesup}}
$$
and
$$
g_{k,s+1} \leq q^{k-1}(q+1)p^{s+1}  \  \  \  \mbox{by Lemma \ref{lemme_genre} (iii).}
$$
So we have
\begin{eqnarray*}
	\frac{g_{k,s+1}}{2X} & = & \frac{g_{k,s+1}}{2(n_0^{k,s}+\frac{D_{k,s}}{2})} \\
				       & \leq & \frac{q^{k-1}(q+1)p^{s+1}}{2(\frac{1}{2}q^{k-1}p^s(q+1)(q-3) + \frac{1}{2}(p-1)p^sq^k)} \\
		    		       & = & \frac{q^{k-1}(q+1)p^sp}{q^{k-1}(q+1)p^s\left(q-3 + (p-1)\frac{q}{q+1}\right) }\\
		    		       & = & \frac{p}{(q-3)+(p-1)\frac{q}{q+1}}
\end{eqnarray*}
Thus, the graph of the function $\Phi$ lies below the line ${y=3\left(1 + \frac{p}{(q-3)+(p-1)\frac{q}{q+1}}\right)x}$. In particular, we get
$$
\Phi(n) \leq 3\left(1 + \frac{p}{(q-3)+(p-1)\frac{q}{q+1}}\right)n.
$$
\item Let $n\geq \frac{1}{2}(p^2+1+{\epsilon (p^2) })$. Otherwise, we already know from Theorems \ref{thm_wdg} and  \ref{thm_shokr} that $\mu_{p^2}(n) \leq 2n$. According to Lemma \ref{lemme_placedegngsr}, there exists a step of the tower $T/\F_{p^2}$ on which we can apply Theorem \ref{thm_arnaud1} with $a=0$. We denote by $H_{k+1}/\F_{p^2}$ the first step of the tower that suits the hypothesis of Theorem~\ref{thm_arnaud1} with $a=0$, i.e. $k$ is an integer such that ${N_{k+1} \geq 2n+2g_{k+1}-1}$ and ${N_k< 2n+2g_k-1}$, where ${N_k:=N_1(H_k/\F_{p^2})}$ and ${g_k:=g(H_k)}$. We denote by $n_0^k$ the biggest integer such that ${N_k\geq 2n_0^k+2g_k-1}$, i.e. \linebreak[4]${n_0^k = \sup \big\{n \in \N \, \vert \, 2n \leq N_k-2g_k+1\big\}}$. To perform multiplication in $\F_{p^{2n}}$, we have the following alternative:
	\begin{enumerate}[(a)]
		\item use the algorithm on the step $H_{k+1}$. In this case, a bound for the bilinear complexity is given by Theorem \ref{thm_arnaud1} applied with $a=0$:
		$$
		\mu_{p^2}(n) \leq 2n+g_{k+1}-1= 2n+g_k-1 +\Delta g_{k,s}.
		$$
		(Recall that $\Delta g_k := g_{k+1} - g_k$)
		\item use the algorithm on the step $H_k$ with an appropriate number of derivative evaluations. Let $a:= 2(n-n_0^k)$ and suppose that $a \leq N_k$. Then ${N_k \geq 2n_0^k+2g_k-1}$ implies that ${N_k +a \geq 2n+2g_k-1}$ so Condition (3) of Theorem \ref{thm_arnaud1} is satisfied. Thus, we can perform $a$ derivative evaluations in the algorithm using the step $H_k$ and we have:
		$$
		\mu_{p^2}(n) \leq 2n+g_k-1+a.
		$$
	\end{enumerate}
	Thus, if $a \leq N_k$ Case (b) gives a better bound as soon as ${a<\Delta g_k}$.
	For $x \in \mathbb{R}^{+}$ such that ${N_{k+1} \geq 2[x]+2g_{k+1}-1}$ and ${N_k < 2[x]+2g_k-1}$, we define the function $\Phi_k(x)$ as follow:
	$$
	\Phi_k(x) = \left\{\begin{array}{ll}
				2x+g_k-1+2(x-n_0^k) & \mbox{if } 2(x- n_0^k) < \Delta g_k\\
				2x+g_{k+1}-1& \mbox{else}.
			      \end{array} \right.
	$$
	Note that when Case (b) gives a better bound, that is to say when  ${2(x-n_0^k) < \Delta g_k}$, then according to Lemma \ref{lemme_deltagsr} we have also 
$${2(x-n_0^k)< N_k}$$
 so we can proceed as in Case (b) since there are enough rational places to use $a=2(x-n_0^k)$ derivative evaluations on.

We define the function $\Phi$ for all ${x\geq0}$ as the minimum of the functions $\Phi_k$ for which $x$ is in the domain of $\Phi_k$. This function is piecewise linear with two kinds of piece: those which have slope $2$ and those which have slope~$4$. Moreover, since the y-intercept of each piece grows with $k$, the graph of the function $\Phi$ lies below any straight line that lies above all the points ${\big(n_0^k+\frac{\Delta g_k}{2}, \Phi(n_0^k+\frac{\Delta g_k}{2})\big)}$, since these are the \textit{vertices} of the graph. Let ${X:=n_0^k+\frac{\Delta g_k}{2}}$, then
	\begin{eqnarray*}
	\Phi(X) & \leq & 2X + g_{k+1} -1 \leq 2\left(1 + \frac{g_{k+1}}{2X}\right)X.
	\end{eqnarray*}
We want to give a bound for $\Phi(X)$ which is independent of $k$.

Lemmas \ref{lemme_genregsr} ii), \ref{lemme_deltagsr} and \ref{lemme_bornesupgsr} give
\begin{eqnarray*}
	\frac{g_{k+1}}{2X} & \leq & \frac{2^{k+2}}{2^{k+1}(p-3)+4+2^{k+1}-2^\frac{k+1}{2}}\\
				     & = &  \frac{2^{k+2}}{2^{k+1}\left((p-3)+1+2^{-k+1}-2^{-\frac{k+1}{2}}\right)}\\
				     & = &  \frac{2}{p-2+2^{-k+1}-2^{-\frac{k+1}{2}}}\\
				     & \leq & \frac{2}{p-\frac{33}{16}}
\end{eqnarray*}
since $-\frac{1}{16}$ is the minimum of the function ${k \mapsto 2^{-k+1}-2^{-\frac{k+1}{2}}}$.\\
Thus, the graph of the function $\Phi$ lies below the line ${y=2\left(1+  \frac{2}{p-\frac{33}{16}}\right)x}$. In particular, we get
$$
\Phi(n) \leq 2\left(1+ \frac{2}{p-\frac{33}{16}}\right)n.
$$
\item Let $n\geq \frac{1}{2}(p+1+{\epsilon (p) })$. Otherwise, we already know from Theorems \ref{thm_wdg} and  \ref{thm_shokr} that $\mu_{p}(n) \leq 2n$. According to Lemma \ref{lemme_placedegngsr2}, there exists a step of the tower $T/\F_p$ on which we can apply Theorem \ref{thm_arnaud2} with $a_1=a_2=0$. We denote by $H_{k+1}/\F_p$ the first step of the tower that suits the hypothesis of Theorem \ref{thm_arnaud2} with $a_1=a_2=0$, i.e. $k$ is an integer such that ${N_{k+1} \geq 2n+2g_{k+1}-1}$ and ${N_k< 2n+2g_k-1}$, where ${N_k:=N_1(H_k/\F_p)+2N_2(H_k/\F_p)}$ and ${g_k:=g(H_k)}$. We denote by $n_0^k$ the biggest integer such that ${N_k\geq 2n_0^k+2g_k-1}$, i.e.\linebreak[4]  ${n_0^k = \sup \big\{n \in \N \, \vert \, 2n \leq N_k-2g_k+1\big\}}$. To perform multiplication in $\F_{p^n}$, we have the following alternative:
	\begin{enumerate}[(a)]
		\item use the algorithm on the step $H_{k+1}$. In this case, a bound for the bilinear complexity is given by Theorem \ref{thm_arnaud2} applied with $a_1=a_2=0$:
		$$
		\mu_q(n) \leq 3n+\frac{3}{2}g_{k+1}= 3n_0^k+\frac{3}{2}g_k+3(n-n_0^k) + \frac{3}{2}\Delta g_k.
		$$
		\item use the algorithm on the step $H_k$ with an appropriate number of derivative evaluations. Let ${a_1+2a_2:= 2(n-n_0^k)}$ and suppose that ${a_1+2a_2 \leq N_k}$. Then ${N_k\geq 2n_0^k+2g_k-1}$ implies that ${N_k +a_1+2a_2 \geq 2n+2g_k-1}$. Thus we can perform $a_1+a_2$ derivative evaluations in the algorithm using the step $H_k$ and we have:
		$$
		\mu_p(n) \leq 3n+\frac{3}{2}g_k+\frac{3}{2}(a_1+2a_2)=3n_0^k+\frac{3}{2}g_k+6(n-n_0^k).
		$$
	\end{enumerate}
	Thus, if $a_1+2a_2 \leq N_{k,s}$ Case (b) gives a better bound as soon as ${n-n_0^{k,s}<\frac{1}{2}\Delta g_{k,s}}$.
	For $x \in \mathbb{R}^{+}$ such that ${N_{k+1} \geq 2[x]+2g_{k+1}-1}$ and ${N_k < 2[x]+2g_k-1}$, we define the function $\Phi_k(x)$ as follow:
	$$
	\Phi_k(x) = \left\{\begin{array}{ll}
				3x+\frac{3}{2}g_k+3(x-n_0^k) & \mbox{if } x- n_0^k < \frac{\Delta g_k}{2}\\
				& \\
				3x+\frac{3}{2}g_{k+1}& \mbox{else}.
			      \end{array} \right.
	$$
	Note that when Case (b) gives a better bound, that is to say when  ${2(x-n_0^k) < \Delta g_k}$, then according to Lemma \ref{lemme_deltagsr} we have also 
$${2(x-n_0^k)< N_k}$$
 so we can proceed as in Case (b) since there are enough places of degree 1 and 2 to use $a_1+a_2=2(x-n_0^k)$ derivative evaluations on.

	We define the function $\Phi$ for all ${x\geq0}$ as the minimum of the functions $\Phi_k$ for which $x$ is in the domain of $\Phi_k$. This function is piecewise linear with two kinds of piece: those which have slope $3$ and those which have slope~$6$. Moreover, since the y-intercept of each piece grows with $k$, the graph of the function $\Phi$ lies below any straight line that lies above all the points ${\big(n_0^k+\frac{\Delta g_k}{2}, \Phi(n_0^k+\frac{\Delta g_k}{2})\big)}$, since these are the \textit{vertices} of the graph. Let ${X:=n_0^k+\frac{\Delta g_k}{2}}$, then
	\begin{eqnarray*}
	\Phi(X) & \leq & 3X + \frac{3}{2}g_{k+1} = 3\left(1 + \frac{g_{k+1}}{2X}\right)X.
	\end{eqnarray*}
We want to give a bound for $\Phi(X)$ which is independent of $k$.

The same reasoning as in (iii) gives 
\begin{equation*}
	\frac{g_{k+1}}{2X}  \leq \frac{2}{p-\frac{33}{16}}
\end{equation*}
Thus, the graph of the function $\Phi$ lies below the line ${y=3\left(1+  \frac{2}{p-\frac{33}{16}}\right)x}$. In particular, we get
$$
\Phi(n) \leq 3\left(1+ \frac{2}{p-\frac{33}{16}}\right)n.
$$

\qed
\end{enumerate}
\end{Proof}

\subsection{New asymptotical upper bounds for $\mu_{q}(n)$}

In this section, we give upper bounds for the asymptotical quantities $m_q$ and $M_q$ which are defined above in Section \ref{mM}. 
First, let us repair the two main mistaken statements (as well as their corollaries) due to I. Shparlinsky, M. Tsfasman and S. Vladut (Theorem 3.1 and Theorem 3.9 in \cite{shtsvl}) in the two following propositions.

\begin{propo}
 Let $q$ be a prime power such that $A(q)>2$. Then
$$m_q \leq 2\left(1+\frac{1}{A(q)-2}\right).$$
\end{propo}

\begin{Proof}\label{newboundmq}
 Let $\left(F_s/\F_q\right)_s$ be a sequence of algebraic function fields defined over $\F_q$.
Let us denote by $g_s$ the genus of $F_s/\F_q$ and by $N_1(s)$ the number of places of degree $1$
of $F_s/\F_q$. Suppose that the sequence $\left(F_s/\F_q\right)_s$ was chosen such that:
\begin{enumerate}
 \item $\lim_{s \rightarrow +\infty}g_s=+\infty$;
 \item $\lim_{s \rightarrow +\infty}\frac{N_1(s)}{g_s}=A(q)$.
\end{enumerate}
Let $\epsilon$ be any real number such that
$0 < \epsilon < \frac{A(q)}{2} -1$.
Let us define the following integer
 $$n_s=\left\lfloor\frac{N_1(s)-2g_s(1+\epsilon)}{2}\right\rfloor.$$
Let us remark that
$$N_1(s)=g_s A(q) + o(g_s),$$
$$\mbox{so }N_1(s)-2(1+\epsilon)g_s=g_s\left(\strut A(q)-2(1+\epsilon)\right)+o(g_s).$$
Then the following holds
\begin{enumerate}
 \item there exists an integer $s_0$ such that for any $s \geq s_0$ the integer $n_s$
is strictly positive;
 \item for any real number $c$ such that $0<c<A(q)-2(1+\epsilon)$ there exists
an integer $s_1$ such that for any integer $s\geq s_1$ the following holds:
$n_s \geq \frac{c}{2}g_s$, hence $n_s$ tends to $+\infty$;
 \item there exists an integer $s_2$ such that for any integer $s\geq s_2$
the following holds: 
$2g_s+1 \leq q^{\frac{n_s-1}{2}}\left(q^{\frac{1}{2}}-1\right)$
and consequently there exists a place of degree $n_s$ (cf. \cite[Corollary 5.2.10 (c) p. 207]{stic} ).
 \item the following inequality holds:
$N_1(s)> 2n_s+2g_s-2$ and consequently, using Theorem 
\ref{theoprinc}  we conclude that $\mu_q(n_s)  \leq  2n_s+g_s-1$.                                                                                   
\end{enumerate}
Consequently, 
$$\frac{\mu_q(n_s)}{n_s} \leq 2+\frac{g_s-1}{n_s},$$
$$m_q \leq 2+ \lim_{s \rightarrow +\infty}\frac{2g_s-2}{N_1(s)-2(1+\epsilon)g_s-2}
\leq 2\left( 1+ \frac{1}{A(q)-2(1+\epsilon)}\right).$$
This inequality is true for any $\epsilon >0$ sufficiently small. Then 
we obtain the result. \qed
\end{Proof}

\begin{coro}\label{coromq1}
Let $q=p^m$ be a prime power such that $q \geq 4$.
Then
$$m_{q^2}\leq 2\left(1+\frac{1}{q-3}\right).$$
\end{coro}

Note that this corollary lightly improves Theorem \ref{chudmq}. Now in the case of arbitrary $q$, we obtain:

\begin{coro}\label{coromq2}
For any $q=p^m>3$,

$$m_{q}\leq 3\left(1+\frac{1}{q-3}\right).$$
\end{coro}

\begin{Proof}
For any $q=p^m>3$, we have $q^2=p^{2m}\geq 16$ and thus Corollary \ref{coromq1} gives 
$m_{q^2}\leq 2\left(1+\frac{1}{q-3}\right)$. Then, by Lemma \ref{lemasyMqmq}, we have 
$$m_q\leq m_{q^2}.\mu_q(2)/2$$
which gives the result since $\mu_q(2)=3$ for any $q$. \qed
\end{Proof}

\vspace{.5em}

Now, we are going to show that for $M_q$ 
the same upper bound as for $m_q$ can be proved though only in the case of $q$ being an even power of a prime. 
However, we are going to prove that in the case of $q$ being an odd power of a prime, the difference between 
the two bounds is very slight.

\begin{propo}\label{newbound}
Let $q=p^m$ be a prime power such that $q \geq 4$. Then
$$M_{q^2}\leq 2\left(1+\frac{1}{q-3}\right).$$
\end{propo}

\begin{Proof}
Let $q=p^m$ be a prime power such that $q\geq4$. Let us consider two cases.
First, we suppose $q=p$. We know that for any real number $\epsilon >0$ and for any sufficiently large real number $x$, 
there exists a prime number $l_k$ such that $x<l_k<(1+\epsilon)x$.
Now, without less of generality let us consider the characteristic $p$ such that $p\neq 11$. 
Then it is known (\cite{tsvl} and \cite{shtsvl})  that the curve $X_k=X_0(11l_k)$, 
where $l_k$ is the $k$-th prime number, has a genus $g_k=l_k$ and satisfies 
$N_1(X_k(\F_{q^2}))\geq (q-1)(g_k+1)$ where $N_1(X_k(\F_{q^2}))$ denotes the number 
of  rational points over $\F_{q^2}$ of the curve $X_k$.
Let us consider a sufficiently large $n$.
There exist two consecutive prime numbers $l_k$ and $l_{k+1}$ such 
that $(p-1)(l_{k+1}+1)> 2n+2l_{k+1}-2$ and $(p-1)(l_k+1)\leq 2n+2l_k-2$. 
Let us consider  the algebraic function field 
$F_{k+1}/\F_{p^2}$ associated to the curve $X_{k+1}$ of genus $l_{k+1}$ defined over $\F_{p^2}$. 
Let $N_i(F_{k}/\F_{p^2})$ be the number 
of places of degree $i$ of $F_{k}/\F_{p^2}$. Then
$N_1(F_{k+1}/\F_{p^2})\geq (p-1)(l_{k+1}+1)> 2n+2l_{k+1}-2$. 
Moreover, it is known that $N_n(F_{k+1}/\F_{p^2})>0$ for any integer $n$ sufficiently large. 
We also know that $l_{k+1}-l_k\leq l_k^{0,535}$ for any integer $k\geq k_0$ where $k_0$ 
can be effectively determined by \cite{baha}. Then there exists a real number $\epsilon>0$ such that 
$l_{k+1}-l_k=\epsilon l_k\leq  l_k^{0,535}$ namely $l_{k+1}\leq(1+\epsilon)l_k$. 
It is sufficient to choose $\epsilon$ such that $\epsilon l_k^{0,465}\leq 1$.
 Consequently, for any integer $n$ sufficiently large, this algebraic function field $F_{k+1}/\F_{p^2}$ satisfies Theorem \ref{theoprinc}, 
and so $\mu_{p^2}(n)\leq 2n+l_{k+1}-1\leq 2n+(1+\epsilon)l_k-1$ with $l_k\leq \frac{2n}{p-3}-\frac{p+1}{p-3}$. 
Thus, as $n\longrightarrow +\infty$ then $l_k\longrightarrow +\infty$ and $\epsilon \longrightarrow 0$, so
we obtain $M_{p^2}\leq 2\left(1+\frac{1}{p-3}\right)$. Note that for $p=11$, Proposition 4.1.20 in \cite{tsvl} 
enables us to obtain $g_k=l_k+O(1)$.

Now, let us study the more difficult case where $q=p^m$ with $m>1$. 
We use the Shimura curves as in \cite{shtsvl}. 
Recall the construction of this good family. Let $L$ be a totally real abelian over 
$\Q$ number field of degree $m$ in which $p$ is inert, thus the residue class field 
${\mathcal O}_L/(p)$ of $p$, where ${\mathcal O}_L$ denotes the ring of integers of $L$, 
is isomorphic to the finite field $\F_{q}$. Let $\wp$ be a prime of $L$ which does 
not divide $p$ and let $B$ be a quaternion algebra for which 
$$B\otimes_{\Q}\R=M_2(\R) \otimes \mathbb{H} \otimes ...\otimes \mathbb{H}$$ 
where $\mathbb{H}$ is the skew field of Hamilton quaternions. Let $B$ be 
also unramified at any finite place if $(m-1)$ is even; let $B$ be also unramified 
outside infinity and $\wp$ if $(m-1)$ is odd. Then, over $L$ one can define the 
Shimura curve by its complex points $X_{\Gamma}(\C)=\Gamma\setminus \mathfrak{h}$, where 
$\mathfrak{h}$ is the Poincar\'e upper half-plane and $\Gamma$ is the group of units 
of a maximal order ${\mathcal O}$ of $B$ with totally positive norm modulo its center. 
Hence,  the considered Shimura curve admits an integral model over $L$ and 
it is well known that its reduction $X_{\Gamma,p}(\F_{p^{2m}})$ modulo $p$ 
is good and is defined over the residue class field ${\mathcal O}_L/(p)$ of $p$, 
which is isomorphic to $\F_q$ since $p$ is inert in $L$. 
Moreover, by \cite{ihar}, the number $N_1(X_{\Gamma,p}(\F_{q^2}))$ of $\F_{q^2}$-points 
of $X_{\Gamma,p}$ is such that $N_1(X_{\Gamma,p}(\F_{q^2}))\geq (q-1)(g+1)$, where $g$ denotes the 
genus of $X_{\Gamma,p}(\F_{q^{2}})$. 
Let now $l$ be a prime which is greater than the 
maximum order of stabilizers $\Gamma_z$, where $z \in \mathfrak{h}$ is a fixed point 
of $\Gamma$ and let $\wp \nmid l$. Let $\Gamma_0(l)_l$ be the following subgroup of $GL_2(\Z_l)$:

$$
\Gamma_0(l)_l=\left \lbrace
\left (
\begin{array}{ll}
 a & b \cr
 c & d
\end{array}
\right )
\in GL_2(\Z_l), c \equiv 0~(mod~l) \right \rbrace .
$$

Suppose that $l$ splits completely in $L$. Then there exists an embedding $F \longrightarrow \Q_l$ where $\Q_l$ denotes 
the usual $l$-adic field, and since $B\otimes_{\Q} \Q_l=M_2(\Q_l)$, we have a natural map: 
$$\phi_l: \Gamma \rightarrow GL_2(\Z_l).$$
Let $\Gamma_l$ be the inverse map of $\Gamma_0(l)_l$ in $\Gamma$ under $\phi_l$. 
Then $\Gamma_l$ is a subgroup of $\Gamma$ of index $l$. We consider the Shimura curve $X_l$ with 
$$X_l(\C)=\Gamma_l\setminus\mathfrak{h}.$$ 
It admits an integral model over $L$ and so can be defined over $L$. Hence, its reduction $X_{l,p}$ modulo $p$ is good 
and it is defined over the residue class field ${\mathcal O}_L/(p)$ of $p$, which is isomorphic to $\F_q$ since $p$ is inert in $L$. 
Moreover the supersingular $\F_p$-points of  $X_{\Gamma,p}$ split completely in 
the natural projection $$\pi_l: X_{l,p} \rightarrow X_{\Gamma,p}.$$ Thus, 
the number of the rational points of $X_{l,p}(\F_{q^2})$ is: 
$$N_1(X_{l,p}(\F_{q^2}))\geq l(q-1)(g+1).$$ Moreover, since $l$ is greater 
than the maximum order of a fixed point of $\Gamma$ on $\mathfrak{h}$, the 
projection $\pi_l$ is unramified and thus by Hurwitz formula, $$g_l=1+l(g-1)$$ 
where $g_l$ is the genus of $X_l$ (and also of $X_{l,p}$).


Note that since the field $L$ is abelian over $\Q$, there exists an integer $N$ such that  field $L$ 
is contained in a cyclotomic extension $\Q(\zeta_N)$ where $\zeta_N$ denotes a primitive root of unity with 
minimal polynomial $\Phi_{N}$. Let us consider the reduction $\Phi_{N,l_k}$ of  $\Phi_{N}$ modulo the prime $l_k$. 
Then, the prime $l_k$ is totally split in the integer ring of $L$ if and only if   the polynomial $\Phi_{N,l_k}$ 
is totally split in $\F_{l_k}=\Z/l_k\Z$ i.e if and only if $\F_{l_k}$ contains the Nth roots of  unity which is 
equivalent to $N\mid l_k-1$. Hence, any prime $l_k$ such that $l_k \equiv 1 \mod N$ is totally split in $\Q(\zeta_N)$ 
and then in $L$. Since $l_k$ runs over primes in an arithmetical progression, the ratio of two consecutive prime numbers 
$l_k \equiv 1 \mod N$ tends to one.

Then for any real number $\epsilon >0$, there exists an integer $k_0$ such that for any integer $k\geq k_0$, 
$l_{k+1}\leq (1+\epsilon)l_k$ where $l_k$ and $l_{k+1}$ are two consecutive prime numbers congruent 
to one modulo $N$. Then there exists an integer $n_{\epsilon}$ such that  
for any integer $n\geq n_{\epsilon}$, the integer $k$, such that the two following inequalities hold $$l_{k+1}(q-1)(g+1)> 2n+2g_{l_{k+1}}-2$$ and $$l_k(q-1)(g+1)\leq 2n+2g_{l_k}-2,$$ 
satisfies $k\geq k_0$ where $g_{l_i}=1+l_i(g-1)$ for any integer $i$.  
Let us consider  the algebraic function field $F_{k}/\F_{q^2}$ defined over the finite field $\F_{q^2}$ 
associated to the Shimura curve $X_{l_{k}}$ of genus $g_{l_{k}}$. Let $N_i(F_{k}/\F_{q^2})$ be the number 
of places of degree $i$ of $F_{k}/\F_{q^2}$.
Then $N_1(F_{k+1})/\F_{q^2})\geq l_{k+1}(q-1)(g+1) > 2n+2g_{l_{k+1}}-2$ where $g$ is the genus 
of the Shimura curve $X_{\Gamma,p}(\F_{q^{2}})$.
Moreover, it is known that there exists an integer $n_0$ such that for any integer $n\geq n_0$,  $N_n(F_{k+1}/\F_{q^2})>0$. 
Consequently, for any integer $n\geq \max(n_{\epsilon},n_0)$ this algebraic function field $F_{k+1}/\F_{q^2}$ 
satisfies Theorem \ref{theoprinc} and so 
$\mu_{q^2}(n)\leq 2n+g_{l_{k+1}}-1\leq 2n+l_{k+1}(g-1)\leq 2n+(1+\epsilon)l_k(g-1)$ with $l_k< \frac{2n}{(q-1)(g+1)-2(g-1)}$. 
Thus, for any real number $\epsilon >0$ and for any $n\geq \max(n_{\epsilon},n_0)$, we obtain $\mu_{q^2}(n)\leq 2n+\frac{2n(1+\epsilon)(g-1)}{(q-1)(g+1)-2(g-1)}$ 
which gives $M_{q^2}\leq2\left(1+\frac{1}{q-3}\right)$. \qed
\end{Proof}

\begin{propo}\label{newbound2}
Let $q=p^m$ be a prime power with odd $m$ such that $q \geq 5$ .
Then
$$M_{q}\leq 3\left(1+\frac{2}{q-3}\right).$$
\end{propo}

\begin{Proof}
It is sufficient to consider the same families of curves that in Proposition \ref{newbound}. 
These families of curves $X_k$ are defined over the residue class field of $p$ which  is isomorphic to $\F_q$.
Hence, we can consider the associated algebraic function fields $F_k/\F_q$ defined over $\F_q$. 
If $q=p$, we have $N_1(F_{k+1}/\F_{p^2})=N_1(F_{k+1}/\F_{p})+2N_2(F_{k+1}/\F_{p})\geq (p-1)(l_{k+1}+1)> 2n+2l_{k+1}-2$ 
since $F_{k+1}/\F_{p^2}=F_{k+1}/\F_{p}\otimes_{\F_p} \F_{p^2}$.  
Then, for any real number $\epsilon >0$ and for any integer $n$ sufficiently large, we have
$\mu_{p}(n)\leq 3n+3g_{l_{k+1}}\leq 3n+3(1+\epsilon)l_k$ by Theorem \ref{theoprinc} 
since $N_n(F_{k+1}/\F_{q^2})>0$. Then, by using the condition $l_k\leq \frac{2n}{p-3}-\frac{p+1}{p-3}$, 
we obtain $M_{p}\leq 3\left(1+\frac{2}{p-3}\right)$.
If $q=p^m$ with odd $m$, we have 
$N_1(F_{k+1}/\F_{q^2})=N_1(F_{k+1}/\F_{q})+2N_2(F_{k+1}/\F_{q})\geq l_{k+1}(q-1)(g+1)> 2n+2g_{l_{k+1}}-2$
since $F_{k+1}/\F_{q^2}=F_{k+1}/\F_{q}\otimes_{\F_q} \F_{q^2}$.
Then, for any real number $\epsilon >0$ and for any integer $n$ sufficiently large as in Proof \ref{newbound}, we have
$\mu_{q}(n)\leq 3n+3g_{l_{k+1}}\leq 3n+3(1+\epsilon)l_k$ by Theorem \ref{theoprinc} 
since $N_n(F_{k+1}/\F_{q^2})>0$. 
Then, by using the condition $l_k< \frac{2n}{(q-1)(g+1)-2(g-1)}$
we obtain $M_{q}\leq 3\left(1+\frac{2}{q-3}\right)$. \qed
\end{Proof}

\begin{propo}\label{newbound3}

$$M_{2}\leq 13.5.$$
\end{propo}

\begin{Proof}
Let $q=p^m=4$. 
We also use the Shimura curves.
Let $L=\Q(\sqrt{d})$ be a totally real quadratic number field such that $d\equiv 1 \mod 8$. 
Then the prime $p=2$ is totally split in $L$ and so the residue class field 
${\mathcal O}_L/(p)$ of $p$, where ${\mathcal O}_L$ denotes the ring of integers of $L$, 
is isomorphic to the finite field $\F_{2}$. Then, let $\wp$ be a prime of $L$ which does 
not divide $p$ and let $B$ be a quaternion algebra for which 
$$B\otimes_{\Q}\R=M_2(\R) \otimes \mathbb{H}$$ 
where $\mathbb{H}$ is the skew field of Hamilton quaternions. 
Let $B$ be also unramified 
outside infinity and $\wp$.
Then, over $L$ one can define the 
Shimura curve by its complex points $X_{\Gamma}(\C)=\Gamma\setminus \mathfrak{h}$, where 
$\mathfrak{h}$ is the Poincar\'e upper half-plane and $\Gamma$ is the group of units 
of a maximal order ${\mathcal O}$ of $B$ with totally positive norm modulo its center. 
Hence,  the considered Shimura curve admits an integral model over $L$ and 
it is well known that its reduction $X_{\Gamma,p}(\F_{p^{2m}})$ modulo $p$ 
is good and is defined over the residue class field ${\mathcal O}_L/(p)$ of $p=2$, 
which is isomorphic to $\F_2$ since $p=2$ is totally split in $L$. 
Moreover, by \cite{ihar}, the number $N_1(X_{\Gamma,p}(\F_{q^2})$ of $\F_{q^2}$-points 
of $X_{\Gamma,p}$ is such that $N_1(X_{\Gamma,p}(\F_{q^2}))\geq (q-1)(g+1)$, where $g$ denotes the 
genus of $X_{\Gamma,p}(\F_{q^{2}})$. 
Let now $l$ be a prime which is greater than the 
maximum order of stabilizers $\Gamma_z$, where $z \in \mathfrak{h}$ is a fixed point 
of $\Gamma$ and let $\wp \nmid l$. Let $\Gamma_0(l)_l$ be the following subgroup of $GL_2(\Z_l)$:

$$
\Gamma_0(l)_l=\left \lbrace
\left (
\begin{array}{ll}
 a & b \cr
 c & d
\end{array}
\right )
\in GL_2(\Z_l), c \equiv 0~(mod~l) \right \rbrace.
$$

Suppose that $l$ splits completely in $L$. Then there exists an embedding $F \longrightarrow \Q_l$ where $\Q_l$ denotes 
the usual $l$-adic field, and since $B\otimes_{\Q} \Q_l=M_2(\Q_l)$, we have a natural map: 
$$\phi_l: \Gamma \rightarrow GL_2(\Z_l).$$
Let $\Gamma_l$ be the inverse map of $\Gamma_0(l)_l$ in $\Gamma$ under $\phi_l$. 
Then $\Gamma_l$ is a subgroup of $\Gamma$ of index $l$. We consider the Shimura curve $X_l$ with 
$$X_l(\C)=\Gamma_l\setminus \mathfrak{h}.$$ 
It admits an integral model over $L$ and so can be defined over $L$. Hence, its reduction $X_{l,p}$ modulo $p=2$ is good 
and it is defined over the residue class field ${\mathcal O}_L/(p)$ of $p=2$, which is isomorphic to $\F_2$ since $p=2$ is totally split in $L$. 
Moreover the supersingular $\F_p$-points of  $X_{\Gamma,p}$ split completely in 
the natural projection $$\pi_l: X_{l,p} \rightarrow X_{\Gamma,p}.$$ Thus, 
the number of the rational points of $X_{l,p}(\F_{q^2})$ is: 
$$N_1(X_{l,p}(\F_{q^2}))\geq l(q-1)(g+1).$$ Moreover, since $l$ is greater 
than the maximum order of a fixed point of $\Gamma$ on $\mathfrak{h}$, the 
projection $\pi_l$ is unramified and thus by Hurwitz formula, $$g_l=1+l(g-1)$$ 
where $g_l$ is the genus of $X_l$ (and also of $X_{l,p}$). 
Note that since the field $L$ is abelian over $\Q$, there exists an integer $N$ such that  field $L$ 
is contained in a cyclotomic extension $\Q(\zeta_N)$ where $\zeta_N$ denotes a primitive root of the unity with 
minimal polynomial $\Phi_{N}$. Let us consider the reduction $\Phi_{N,l_k}$ of  $\Phi_{N}$ modulo the prime $l_k$. 
Then, the prime $l_k$ is totally split in the integer ring of $L$ if and only if   the polynomial $\Phi_{N,l_k}$ 
is totally split in $\F_{l_k}=\Z/l_k\Z$ i.e if and only if $\F_{l_k}$ contains the Nth roots of the unity which is 
equivalent to $N\mid l_k-1$. Hence, any prime $l_k$ such that $l_k \equiv 1 \mod N$ is totally split in $\Q(\zeta_N)$ 
and then in $L$. Since $l_k$ runs over primes in an arithmetical progression,  the ratio of two consecutive prime numbers 
$l_k \equiv 1 \mod N$ tends to one.
Then for any real number $\epsilon >0$, there exists an integer $k_0$ such that for any integer $k\geq k_0$, 
$l_{k+1}\leq (1+\epsilon)l_k$ where $l_k$ and $l_{k+1}$ are two consecutive prime numbers congruent 
to one modulo $N$. Then there exists an integer $n_{\epsilon}$ such that  for any integer $n\geq n_{\epsilon}$, 
the integer $k$, such that the two following inequalities hold $$l_{k+1}(q-1)(g+1)> 2n+2g_{l_{k+1}}+6$$ 
and $$l_k(q-1)(g+1)\leq 2n+2g_{l_k}+6,$$ satisfies $k_\geq k_0$ where $g_{l_i}=1+l_i(g-1)$ for any integer $i$.

Let us consider  the algebraic function field $F_{k}/\F_{2}$ defined over the finite field $\F_{2}$ 
associated to the Shimura curve $X_{l_{k}}$ of genus $g_{l_{k}}$. Let $N_i(F_{k}/\F_{t})$ be the number 
of places of degree $i$ of $F_{k}/\F_{t}$ where $t$ is a prime power.
Then, since $F_{k+1}/\F_{q^2}=F_{k+1}/\F_{2}\otimes_{\F_2} \F_{q^2}$ for $q=4$, we have 
$N_1(F_{k+1}/\F_{q^2})=N_1(F_{k+1}/\F_{2})+2N_2(F_{k+1}/\F_{2})+4N_4(F_{k+1}/\F_{2})\geq l_{k+1}(q-1)(g+1)> 2n+2g_{l_{k+1}}+6$ 
where $g$ is the genus of the Shimura curve $X_{\Gamma,p}(\F_{q^{2}})$. 
Moreover, it is known that there exists an integer $n_0$ such that for any integer $n\geq n_0$,  $N_n(F_{k+1}/\F_{q^2})>0$. 
Consequently, for any integer $n\geq \max(n_{\epsilon},n_0)$ this algebraic function field $F_{k+1}/\F_{2}$ 
satisfies Theorem 3.2 in \cite{bapi} and so 
$\mu_{2}(n)\leq \frac{9}{2}(n+g_{l_{k+1}}+5)\leq \frac{9}{2}(n+l_{k+1}(g-1)+6)\leq \frac{9}{2}(n+(1+\epsilon)l_k(g-1))+27$ 
with $l_k< \frac{2n+8}{(q-1)(g+1)-2(g-1)}$. 
Thus, for any real number $\epsilon >0$ and for any $n\geq \max(n_{\epsilon},n_0)$, we obtain 
$\mu_{2}(n)\leq \frac{9}{2}(n+2n\frac{(1+\epsilon)}{q-3}+\frac{8}{q-3})+27\leq \frac{9}{2}(1+2(1+\epsilon))n+63$  
which gives $M_{2}\leq 13,5$. \qed

\end{Proof}


\end{document}